 \documentclass[reqno,11pt]{amsart}
 \usepackage{amsmath, amsthm, a4, latexsym, amssymb}

\def\@rmrk#1#2{\refstepcounter
    {#1}\@ifnextchar[{\@yrmrk{#1}{#2}}{\@xrmrk{#1}{#2}}}

\makeatletter\@addtoreset{equation}{section}\makeatother

 \newfont{\bfit}{cmbxti10 scaled 1200}

 \renewcommand{\d}{d}
 
 \newcommand{\En}{{\rm E} }
 \newcommand{\Sn}{{\mathbb S}_n}
 \newcommand{\eps}{\varepsilon}
 \newcommand{\supp}{{\rm supp}}

 \newcommand{\R}{\mathbb{R}}
 \newcommand{\N}{\mathbb{N}}
 
 \newcommand{\Sym}{\mathfrak{S}}
 \newcommand{\prob}{\mathbb{P}}
 
 \newcommand{\Leb}{{\rm Leb\,}}
 
 \newcommand{\me}{\mathbb{E}}
 \newcommand{\E}{\mathbb{E}}
 \renewcommand{\P}{\mathbb{P}}
 \newcommand{\1}{{\sf 1}}
 \newcommand{\tmu}{\widetilde{\mu}}

 \newcommand{\och}{\mathfrak{H}}
 \newcommand{\A}{{\mathfrak A}}

 \newcommand{\skric}{{\mathcal C}}
 \newcommand{\skrid}{{\mathcal D}}
 
 \newcommand{\skrif}{{\mathcal F}}
 \newcommand{\skrig}{{\mathcal G}}

 \newcommand{\skrim}{{\mathcal M}}
 
 \newcommand{\skrip}{{\mathcal P}}

\newcommand{\heap}[2]{\genfrac{}{}{0pt}{}{#1}{#2}}
\newcommand{\sfrac}[2]{\mbox{$\frac{#1}{#2}$}}
\newcommand{\ssup}[1] {{\scriptscriptstyle{({#1}})}}

\newenvironment{Proof}[1]
{\vskip0.1cm\noindent{\bf #1}}{\vspace{0.15cm}}
\renewcommand{\subsection}{\secdef \subsct\sbsect}
\newcommand{\subsct}[2][default]{\refstepcounter{subsection}
\vspace{0.15cm}
{\flushleft\bf \arabic{section}.\arabic{subsection}~\bf #1  }
\nopagebreak\nopagebreak}
\newcommand{\sbsect}[1]{\vspace{0.1cm}\noindent
{\bf #1}\vspace{0.1cm}}

{\nopagebreak {\hfill\rule{2mm}{2mm}}\\ }

\newtheorem{theorem}{Theorem}[section]
\newtheorem{lemma}[theorem]{Lemma}
\newtheorem{cor}[theorem]{Corollary}
\newtheorem{prop}[theorem]{Proposition}

\newtheoremstyle{thm}{1.5ex}{1.5ex}{\itshape\rmfamily}{}
{\bfseries\rmfamily}{}{2ex}{}

\newtheoremstyle{rem}{1.3ex}{1.3ex}{\rmfamily}{}
{\itshape\rmfamily}{}{1.5ex}{}
\theoremstyle{rem}
\newtheorem{remark}{{\slshape\sffamily Remark}}[]

\refstepcounter{subsubsection}

\def\thebibliography#1{\section*{References}
  \list%
  {\arabic{enumi}.}
    {\settowidth\labelwidth{[#1]}\leftmargin\labelwidth
    \advance\leftmargin\labelsep
    \parsep0pt\itemsep0pt
    \usecounter{enumi}}
    \def\newblock{\hskip .11em plus .33em minus .07em}
    \sloppy                   
    \sfcode`\.=1000\relax}

\begin{document}
\title[Exponential moments for intersection local times]
{\large Brownian intersection local times:\\ Exponential moments and law
of large masses}
\author[Wolfgang K\"onig and Peter M\"orters]{}
\maketitle
\thispagestyle{empty}
\vspace{-0.5cm}

\centerline{\sc By Wolfgang K\"onig
and Peter M\"orters}
\renewcommand{\thefootnote}{}
\footnote{\textit{AMS Subject
Classification:} 60J65, 60J55, 60F10.}
\footnote{\textit{Keywords:} Intersection of Brownian paths, 
intersection local time, exponential moment, Feynman-Kac formula.}

\vspace{-0.5cm}
\centerline{\textit{Technische Universit\"at Berlin, and University of Bath}}
\vspace{0.2cm}

\begin{quote}{\small {\bf Abstract: }
Consider $p$ independent Brownian motions in $\R^d$, each 
running up to its first exit time from an open domain $B$,
and their intersection local time $\ell$ as a measure on $B$. 
We give a sharp criterion for the finiteness of exponential moments, 
$$\me\Big[\exp\Big(\sum_{i=1}^n 
\langle \varphi_i, \ell \rangle^{1/p}\Big) \Big],$$
where $\varphi_1, \dots,  \varphi_n$ are nonnegative, bounded functions with compact 
support in $B$. We also derive a law of large numbers for 
intersection local time conditioned to have large total mass.}\end{quote}


\section{Introduction and statement of results}\label{intro}

\subsection{Introduction}\label{aims}

Much of our knowledge about Brownian occupation times goes back to the
celebrated \emph{Feynman-Kac formula}, see e.g.\ \cite{Fe48, Ka49} for
original sources and \cite{FP99} for an excellent survey. The formula is
also at the heart of results relating Brownian motion and a vast number
of differential equations, see \cite{BS02} for an impressive account.  
To formulate one of many versions of the Feynman-Kac formula, let $B\subset\R^d$ 
be a bounded open domain and suppose that $W$ is a Brownian motion started in 
$x\in B$. 
Let $q\colon B\to[0,\infty)$ be a Borel measurable function and $T$ the first 
exit time of the Brownian motion $W$ from $B$. Then the function 
$f\colon B\to[0,\infty]$ given by
\begin{equation}\label{FKform}
f(x)=\me_x\Big[ \exp \int_0^T q(W(s)) \, ds \Big],
\mbox{ for } x\in B,
\end{equation}
is the minimal positive solution of the equation
\begin{equation}\label{Greeneq}
f(x)=1+ \int G(x,y) f(y) q(y) \, dy,
\end{equation}
where $G$ is the Green function for Brownian motion killed on the
boundary of $B$, see e.g.\ \cite[(8)]{FP99}. Hence, the exponential moments in
\eqref{FKform} are finite if and only if there exists a {\em finite\/} positive
solution to \eqref{Greeneq}.

An \emph{explicit} criterion for the finiteness of the exponential moments 
in \eqref{FKform}, given in terms of a variational formula, is due to  
Pinsky \cite{Pi86}. Assuming that $B$ is a bounded, open domain with a 
\mbox{$C^2$-boundary}, for any continuous 
function $q\colon \overline{B}\to \R$, define
$$l_{q,B}=\sup\Big\{\int q(x)\psi(x)^2\,dx- \sfrac 12 \|\nabla \psi\|^2_2 \, 
: \, \psi\in H_0^1(B),\int \psi(x)^2\,dx  =1 \Big\}.$$
Pinsky shows that
\begin{equation}\label{pinsky}
\me_x\Big[ \exp \int_0^T q(W(s)) \, ds \Big] \quad
\left\{ \begin{array}{l}
<\infty \mbox{ for all $x\in B$ if } l_{q,B}<0,\\
=\infty \mbox{ for all $x\in B$ if } l_{q,B}>0.\\
\end{array}
\right.
\end{equation}
His main tool is the \emph{Donsker-Varadhan large deviation theory} for the
occupation measure of Brownian motion.

Whereas we have an almost complete understanding of occupation times 
of \emph{one} Brownian path, we know much less about the intersection of 
\emph{several} independent Brownian paths, where the role of occupation 
times is played by \emph{intersection local times}. A major reason for this 
lack of understanding is that we do not know the natural analogues of the two 
crucial tools in our understanding of occupation measures, 
the {Feynman-Kac} formula and the {Donsker-Varadhan} 
large deviation theory. However, there is every indication that the relation 
of intersection local times and differential equations is as rich and 
exciting as in the case of a single Brownian motion. Maybe even more so, 
as the equations appearing in this context turn out to be \emph{nonlinear} 
and the analytical theory of these equations knows more open problems than answers.

In this paper, as a first step into this new territory, we investigate the 
existence of exponential moments for intersection local times of $p$ Brownian
motions in $\R^d$. In analogy to \eqref{pinsky} we give a finiteness criterion 
for exponential moments of integrals of intersection local time against a large 
class of test functions. In the absence of the two crucial tools mentioned
before, our arguments rely heavily on combinatorial and analytical methods.
As a consequence of our approach, we are able to prove a law of large numbers 
(or, more accurately, large masses) which relates
the asymptotic shape of intersection local time to a natural nonlinear partial
differential  equation.

\subsection{Brownian intersection local times}\label{ass}

Suppose that $B\subset\R^d$, with $d\geq 1$, is an open set. 
The set $B$ is assumed to be bounded if $d \le 2$, except that we allow $B=\R^d$, 
if $d\ge 3$. Let $p\ge 1$ be an integer, $x=(x_1,\ldots,x_p)\in B^p$, and
assume that a family of $p$ independent Brownian motions
$$(W^1(t)\colon t\in[0,\infty)),\dots, (W^p(t)\colon t\in[0,\infty))$$
in $\R^d$ with $W^1(0)=x_1, \ldots, W^p(0)=x_p$ are realized on a probability space
$(\Omega, \skrif, \prob_x)$. Denote the corresponding expectation by 
$\me_x$. Each motion is killed at the first exit time $T^i=\inf\{t>0\colon W^{i}(t)
\not\in B\}$ from $B$ if this time is finite. 
By classical results  of Dvoretzky, Erd\H{o}s, Kakutani and Taylor,
$p(d-2)<d$ is equivalent to the fact that the paths of the $p$ motions 
have a positive probability of intersecting in a point other than
their starting point. In this case there exists a  locally finite measure, 
the {\it (projected) intersection local time measure\/} $\ell$, which can 
symbolically  be described by the formula
\begin{equation}\label{localtimedef}
\ell(A)=\int_A dy \, \prod_{i=1}^p \int_0^{T^i}\!\!\d s \, \delta_y
\big(W^{i}(s)\big),
\mbox{ for $A\subset\R^d$ Borel.}
\end{equation}
Heuristically, $\ell(A)$ measures the amount of intersection of the $p$
Brownian paths in the set $A$ before they are killed. The measure $\ell$ is 
a random element of the space $\skrim(B)$ of nonnegative, locally finite 
measures on~$B$, which is equipped with the vague topology. $\ell$ is nontrivial
with positive probability and, if $A$ is a bounded set, 
$\ell(A)$ is almost surely finite. If $A$ is unbounded, then $\ell(A)$ 
may be equal to $\infty$ with positive probability. 

We \emph{always} assume that $p(d-2)<d$. This includes the following cases,
\begin{itemize}
\item $p=1$, $d$ arbitrary. In this case $\ell$ degenerates to the occupation 
measure of a single Brownian path,
$$\ell(A)=\int_0^{T^1}\!\!ds \, \1_A\big(W^{1}(s)\big),
\mbox{ for $A\subset\R^d$ Borel.}$$
Our main result, Theorem~\ref{moments}, 
is essentially contained in \cite{Pi86}, see Remark~\ref{Picompare} for a 
comparison. Theorem~\ref{LLM} seems to be new even in this case. 
\item $d=1$, $p\ge 2$ arbitrary. In this case the symbolic formula 
\eqref{localtimedef} makes sense using local time. Indeed, if 
$(L^i(x)\,:\, x\in \R)$ is the family of local times of the stopped Brownian
motion $W^i$, i.e., the continuous density of the occupation measure
$\int_0^{T^i} ds\, \delta_{W^i(s)}$, we define
$$\ell(A)=\int_A dy \, \prod_{i=1}^p L^i(y),
\mbox{ for $A\subset\R$ Borel.}$$
\item $d=2$, $p\ge 2$ arbitrary, and $d=3$, $p=2$. In these most interesting
cases, the local times do not exist and substantial work is needed to turn 
\eqref{localtimedef} into a rigorous definition.  
See Section~2 of \cite{KM02} for a short survey on three rigorous 
constructions of $\ell$ in these cases.
\end{itemize}
We would like to mention (see \cite{LG87, LG89}) that if $d\ge 2$, almost surely, 
$\ell$ is equal to a Hausdorff measure on the set 
$S:=W^1[0,T^1)\cap\ldots\cap W^p[0,T^p)$ 
with some deterministic gauge function. This fact underlines that $\ell$ is the 
natural measure on the intersection of the paths.

\subsection{The main result}\label{sec-main}

Let $\varphi\colon B\to[0,\infty)$ be bounded with compact support in $B$. If
$\mu$ is a measure on $B$ we write $\langle \varphi, \mu\rangle=\int \varphi 
\, d\mu$. Suppose now that $\varphi$ is positive on a set of positive 
Lebesgue measure, then it turns out that $\me_x[ \exp( \langle \varphi, 
\ell \rangle )]=\infty$. However, it is a subtle question whether the
\emph{stretched} exponential moments of the form
$$\me_x\Big[ \exp\big( \langle \varphi, \ell \rangle^{1/p} \big)\Big]$$
are finite or not. Our first main result is a sharp criterion for this. In fact,
the nonlinearity due to the $p$th root in the exponent makes it natural
to ask a more general question, namely when for a finite family
$(\varphi_1, \ldots, \varphi_n)$ of bounded nonnegative functions the moments
$$\me_x\Big[\exp\Big(\sum_{i=1}^n \langle \varphi_i, \ell \rangle^{1/p}\Big) \Big]$$
are infinite or not. To formulate our answer denote by
\begin{equation}\label{Ddef}
\skrid(B)=\begin{cases} H_0^1(B)&\mbox{if }B \mbox{ is bounded},\\
D^1(\R^d)& \mbox{if }B=\R^d,
\end{cases}
\end{equation}
the classical Sobolev space $H_0^1(B)$ with zero boundary condition if
$B$ is bounded, and, in the case that $B=\R^d$,  
the set $D^1(\R^d)$ of functions in $L^1_{{\rm loc}}(\R^d)$ vanishing at infinity
and having a distributional gradient in $L^2(\R^d)$. In Section~\ref{space}
we recall some properties of $\skrid(B)$.

\begin{theorem}[Exponential moments]\label{moments}
Let $\phi=(\phi_1, \ldots, \phi_n)$ be a family of bounded nonnegative Borel 
measurable functions with compact support in $B$, and let
\begin{equation}\label{Theta}
\Theta(\phi)=\Theta(\phi_1, \ldots, \phi_n) =
\inf\Big\{ \frac p2 \|\nabla \psi\|^2_2 \, : \, \psi\in\skrid(B), \,
\sum_{i=1}^n \| \phi_i \psi \|_{2p}^2 =1 \Big\}.
\end{equation}
Then
\begin{equation}\label{criterion}
\me_x\Big[\exp\Big(\sum_{i=1}^n 
\langle \phi_i^{2p}, \ell \rangle^{1/p}\Big) \Big] 
\left\{
\begin{array}{l}
< \infty \mbox{ for all $x\in B^p$ if } \Theta(\phi)>1, \\
= \infty \mbox{ for all $x\in B^p$ if } \Theta(\phi)<1. \\
\end{array} \right.
\end{equation}
Indeed, we even have 
\begin{equation}\label{tailbehaviour}
\lim_{a\uparrow\infty} \frac{1}{a} \log
\prob_x\Big\{ \sum_{i=1}^n \langle \phi_i^{2p}, \ell \rangle^{1/p} >a\Big\} 
=-\Theta(\phi).
\end{equation}
\end{theorem}

A partial result in the direction of Theorem~\ref{moments} was obtained 
in~\cite{KM02}. {F}rom Theorem~\ref{moments} we can infer a finiteness criterion
for intersection local times in a form analogous to Pinsky's result for single 
Brownian motion in \eqref{pinsky}.

\begin{cor} Let $B\subset\R^d$ be a bounded, open domain with
$C^1$-boundary, and let $\phi=(\phi_1, \ldots, \phi_n)$ be a family of bounded
Borel measurable
functions \mbox{$\phi_i\colon\overline{B}\to[0,\infty)$.} Let
\begin{equation}\label{varPinsky}
l^p_{\phi,B}=\sup\bigg\{\sum_{i=1}^n \| \phi_i \psi \|^2_{2p}- 
\frac p2 \|\nabla \psi \|^2_2\, : \, 
\psi\in\skrid(B), \, \|\psi\|_{2p}=1\bigg\}. 
\end{equation}
Then
\begin{equation}\label{criterionalapinski}
\me_x\Big[\exp\Big(\sum_{i=1}^n 
\langle \phi_i^{2p}, \ell \rangle^{1/p}\Big) \Big] 
\left\{
\begin{array}{l}
< \infty \mbox{ for all $x\in B^p$ if } l^p_{\phi,B}<0, \\
= \infty \mbox{ for all $x\in B^p$ if } l^p_{\phi,B}>0. \\
\end{array} \right.
\end{equation}
\end{cor}

\begin{Proof}{Proof.}
For $\eps>0$, denote by $B[\eps]$ the open $\eps$-neighbourhood of $B$.
We first show that
\begin{equation}
\label{ellcont}
\limsup_{\eps\downarrow 0}  \, l^p_{\phi,B[\eps]} \le l^p_{\phi,B}.
\end{equation}
Indeed, for $k\in\N$ let $\psi_k\in H_0^1(B[1/k])$ be an approximate minimizer
in \eqref{varPinsky} i.e., $\|\psi_k\|_{2p}=1$ and $\sum_{i=1}^n \| \phi_i \psi_k \|^2_{2p}- 
\frac p2 \|\nabla \psi_k \|^2_2\geq l^p_{\phi,B[1/k]}-1/k$. Since the first term is bounded in $k\in\N$,
it is clear that $(\|\nabla \psi_k \|^2)_{k\in\N}$ is bounded. 
By Lemma~\ref{subsequences}
we may assume that $\psi_k $ converges, as $k\to\infty$, to 
some $\psi\in H_0^1(B[1])$
in $L^{2p}$-norm such that $\nabla \psi_k$ converges weakly to $\nabla \psi$.
Since $\supp(\psi_k)\subset B[1/k]$ for any $k\in\N$, we may assume that 
$\psi\in H^1(\R^d)$
with $\psi=0$ outside $B$. According to Lemma~\ref{HaNull}, the restriction of
$\psi$ to $B$ lies in $H_0^1(B)$. By lower semicontinuity of $\|\cdot\|_2$, we have $\|\nabla\psi\|_2
\leq \liminf_{k\uparrow \infty}\|\nabla\psi_k\|_2$. By $L^{2p}$-convergence we have
$\lim_{k\uparrow \infty}\sum_{i=1}^n \| \phi_i \psi_k \|^2_{2p}=\sum_{i=1}^n \| 
\phi_i \psi \|^2_{2p}$. Hence, 
$$
\limsup_{\eps\downarrow0}\, l^p_{\phi,B[\eps]} \leq\sum_{i=1}^n \| 
\phi_i \psi \|^2_{2p}-\frac p2 \|\nabla\psi\|_2^2\leq  l^p_{\phi,B},
$$
and this finishes the proof of \eqref{ellcont}.

Now assume that $l^p_{\phi,B}<0$.
Because of \eqref{ellcont} one can fix $\eps>0$ such that $l^p_{\phi,B[\eps]}<0$.
We now work in the domain $B[\eps]$ and exploit that the supports of 
$\phi_1,\ldots,\phi_n$ are strictly inside $B[\eps]$. For any $\eta>0$,
\begin{equation}\label{corcalc1}
\begin{aligned}
0<-l^p_{\phi,B[\eps]}&=\inf_{\psi\in\skrid(B[\eps])\setminus\{0\}}\frac{\frac p2
 \|\nabla \psi \|^2_2-\sum_{i=1}^n \| \phi_i \psi \|^2_{2p}}{\|\psi\|_{2p}^2}\\
&\leq \inf\Bigl\{\frac{\frac p2 \|\nabla \psi \|^2_2-1}{\|\psi\|_{2p}^2}\,:\,
\psi\in\skrid(B[\eps]),\,
\sum_{i=1}^n \| \phi_i \psi \|^2_{2p}=1,\|\psi\|_{2p}^2\geq \eta\Bigr\}\\
&\leq \frac 1\eta\Bigl(\inf\Bigl\{\frac p2 \|\nabla \psi \|^2_2\,:\,
\psi\in\skrid(B[\eps]),
\,\sum_{i=1}^n \| \phi_i \psi \|^2_{2p}=1,\|\psi\|_{2p}^2\geq \eta\Bigr\}-1\Bigr).
\end{aligned}
\end{equation}
Since there is a non-trivial minimiser $\psi$ for the variational formula 
\eqref{Theta}, we can choose $\eta>0$ so small that the variational 
formula on the right hand side of \eqref{corcalc1} is equal to $\Theta(\phi)$. 
Hence, $\Theta(\phi)>1$, and by Theorem~\ref{moments} we infer that
$$\me_x\Big[\exp\Big(\sum_{i=1}^n 
\langle \phi_i^{2p}, \ell \rangle^{1/p}\Big) \Big]<\infty$$
for all $x\in B[\eps]^p$ and the intersection local time $\ell$ of the
Brownian motion killed upon leaving $B[\eps]$. The desired result for the
original domain follows by monotonicity.

Now assume that $l^p_{\phi,B}>0$. Then there is a $\psi\in\skrid(B)$ satisfying 
$\|\psi\|_{2p}=1$ and $$\frac p2
 \|\nabla \psi \|^2_2<\sum_{i=1}^n \| \phi_i \psi \|^2_{2p}.$$ 
There exists a domain $U$ with $\overline{U}\subset B$ such that
$$\frac p2  \|\nabla \psi \|^2_2<\sum_{i=1}^n \| \phi_i\1_{U} \psi \|^2_{2p}.$$ 
Multiplying $\psi$ with an appropriate positive constant, we obtain a new 
$\psi\in\skrid(B)$ satisfying $\sum_{i=1}^n \| \phi_i\1_{U} \psi \|^2_{2p}=1$ and 
$\frac p2 \|\nabla\psi \|^2_2<1$, which implies that $\Theta(\phi\1_{U})<1$.
The result follows from Theorem~\ref{moments} and monotonicity.
\end{Proof}
\qed

\begin{remark}[Comparison with the result of \cite{Pi86}.]\label{Picompare}
Looking at the special case $p=1$, and (without loss of 
generality)~$n=1$, Pinsky shows \eqref{criterionalapinski} for any 
continuous $\phi_1\colon\overline{B}\to\R$ 
in the case that $B$ has a $C^2$-boundary. 
His result is based on the formula 
\begin{equation}\label{DV}
\lim_{t\uparrow\infty} \frac{1}{t} \, \log 
\me_x\Big[ \exp\Big\{\int_0^t \phi_1(W(s))\, ds\Big\},
\, T^1>t \Big] = l^1_{\phi_1,B}, \mbox{ for all }x\in B,
\end{equation}
which follows from the {Donsker-Varadhan theory} and an analysis of the
regularity of the variational problem defining $l^1_{\phi,B}$. 
Note that our approach requires weaker regularity assumptions, but
is not suitable to deal with functions of changing sign.
Of course, the main point of our investigation is the generalisation to the case 
$p>1$, where \eqref{DV} is not available. \hfill{$\Diamond$}
\end{remark}

In our proof of Theorem~\ref{moments} it is essential to show the 
\emph{existence} of minimisers in \eqref{Theta} and characterise them
by differential equations.

\begin{prop}[Analysis of $\Theta(\phi)$]\label{analtheta}
Let $\phi=(\phi_1, \ldots, \phi_n)$ be a family of nonnegative bounded  measurable 
functions with compact support in $B$. Then the infimum in \eqref{Theta}
is attained. Every minimiser $\psi \in\skrid(B)$ satisfies the equation
\begin{equation}\label{nonlinear}
-\frac p2 \Delta\psi = \Theta(\phi)\,
\psi^{2p-1} \sum_{i=1}^n \|\phi_i \psi\|_{2p}^{2-2p}
\, \phi_i^{2p}.
\end{equation}
\end{prop}

We do not know in general whether the solution of \eqref{nonlinear} is 
\emph{unique}, even in the simpler case $n=1$, where the equation simplifies to
the nonlinear eigenvalue equation 
$$-\frac p2 \Delta\psi = \Theta(\phi)\, \psi^{2p-1} \phi^{2p}.$$
A natural question in this context is whether the minimisers in the variational 
problem \eqref{Theta} allow a probabilistic interpretation. Our
second main result provides an interpretation as the asymptotic ``shape'' of
the intersection local time when the total mass in a given set is large.

For the formulation of this result let $U$ be an open, bounded Lebesgue
continuity set whose closure is contained in $B$, and define a (random) 
probability measure $L$ on $U$ as the normalized intersection local times 
$L=\ell/\ell(U)$ on $U$.
Let $\sf{d}$ denote a metric on the space $\skrim_1(U)$ of probability measures
on $U$, that induces the weak topology. Under the conditional law 
$\P\{\,\cdot\mid\ell(U)>a\}$, as $a\uparrow \infty$,
the measure $L$ satisfies the following law of large numbers.

\begin{theorem}[Law of Large Masses] \label{LLM}
Denote by $\mathfrak{M}\subset\skrim_1(U)$ the set of measures 
$\psi^{2p}(x)\, dx$ on $U$ with $\psi$ a minimiser in the variational 
formula for $\Theta(\1_U)$. Then, for any $x\in B^p$,
$$ \lim_{a\uparrow\infty}\P_x\Big\{{\sf{d}}\big( L, \, 
\mathfrak{M}\big)\ge \eps\,\Big| \, \ell(U)>a\Big\}=0,
\qquad\mbox{ for any $\eps>0$.}$$
The convergence is exponential with speed $a^{1/p}$.
\end{theorem}

This result was announced, without proof, in \cite{KM02}.
By Proposition~\ref{analtheta}, the densities $\psi^{2p}$ of the
measures in $\mathfrak{M}$ satisfy the nonlinear eigenvalue equation
$$- \frac p2 \Delta \psi = \Theta(\1_U)\, \psi^{2p-1} \1_U.$$
Whether the solution to this equation is \emph{unique}, and also whether
$\mathfrak{M}$ is a singleton, seems to be an {open} problem. It is
a further {open} problem to determine the precise rate of convergence 
in Theorem~\ref{LLM}.

\begin{remark}[Comparison with a result of \cite{Ch03}.]\label{Chcompare}
In a recent preprint Xia Chen~\cite{Ch03} looks at the intersection 
local time $\ell_1$ for $p>1$ Brownian paths in $\R^d$ each running
up to time one. He finds a sharp criterion for finiteness of the 
total intersection local time. For $q=\sfrac{d}{2} (p-1)$ and 
$$\gamma(p,d)= q\big(\sfrac{p}{p-q}\big)^{1- \frac pq}
\sup\Big\{ \|\psi\|^2_{2p} - \sfrac{1}{2} \| \nabla \psi \|^2_2 \, : 
\, \psi\in \skrid(\R^d),  \|\psi\|_2=1\Big\}^{1-\frac{p}{q}},$$
he shows that
\begin{equation}\label{chencriterion}
\me_x\Big[\exp\big( \gamma \ell_1(\R^d)^{2/d(p-1)}\big) \Big] 
\left\{ \begin{array}{l}
< \infty \mbox{ for all $x\in B^p$ if } \gamma<\gamma(p,d), \\
= \infty \mbox{ for all $x\in B^p$ if } \gamma>\gamma(p,d). \\
\end{array} \right.
\end{equation}
Because of the \emph{fixed} time horizon this problem is very different 
from the one we are looking at --- note also the completely different 
scaling behaviour. Still it is interesting to compare the techniques 
of proof. Chen's method, again quite different from ours, is based on asymptotics 
for the Feynman-Kac formula for occupation times and approximation of intersection
local time by occupation time. Understanding the relationship between these 
results and methods would probably lead to significant progress in the research 
programme set out in the introduction.\hfill{$\Diamond$}
\end{remark}

\section{Overview and setup of the proof}

Our main results follow from an analysis of the large-$k$ asymptotics of
the $k$th moments of the random variables 
$$\sum_{i=1}^n \langle \phi^{2p}_i, \ell \rangle^{1/p}.$$
In Section~\ref{first} we fix some notation about entropy and relative entropy.
In Section~\ref{subsec-momasy} we derive the moment asymptotics in terms of a 
variational formula involving relative entropies. We also identify this formula 
in terms of $\Theta(\phi)$ defined in \eqref{Theta}. In Section~\ref{sec-finish}
we complete the proofs.

\subsection{Entropy and relative entropy}\label{first}

For any probability measure $\mu$ and any measure $\tmu$ on the
same measure space $X$ the \emph{relative entropy} or \emph{Kullback-Leibler 
distance} of $\mu$ with respect to $\tmu$ is defined as
\begin{equation}\label{entropy}
H(\mu \mid \tmu )= \int \mu(\d x) \, \log\frac{\mu(\d x)}{\tmu(\d x)},
\end{equation}
which is to be interpreted as infinity if $\mu\not\ll\tmu$. If $\tmu=f\, dx$
then we often write $H(\mu \mid f)$ instead of $H(\mu\mid\tmu)$. By Jensen's 
inequality we always have $H(\mu\,|\,\tmu)\ge -\log \tmu(X)$ and equality 
holds if and only if $\mu=\tmu/\tmu(X)$. 
More specifically, for any $\mu\in\skrim_1(B)$ we define $I(\mu)=
H(\mu\mid\Leb)$, the relative entropy of $\mu$ with respect to
the Lebesgue measure $\Leb$ on~$B$. For $\lambda\in\mathbb{S}_n$,
where $\mathbb{S}_n$ denotes the $n$-dimensional unit simplex,
we write $I(\lambda)=\sum_{i=1}^n \lambda_i \log \lambda_i$ for the 
relative entropy of $\lambda$ with respect to the counting measure. 
In both cases the functional $I$ is convex and lower semicontinuous. 
We denote by
\begin{equation}\label{M2def}
\skrim^*_1(B)= \Big\{ \nu\in\skrim_1(B^2) \, : \, \nu(A\times B)=\nu( B\times A)
\mbox{ for all Borel sets }A\subset B  \Big\}
\end{equation}
the set of probability measures $\nu$ on $B^2$ with equal marginals
$\nu_1(A)=\nu(A\times  B)$ and $\nu_2(A) =\nu( B \times A)$, we denote 
the marginal by $\bar{\nu}=\nu_1=\nu_2$. For $\nu\in\skrim_1(B^2)$ we define
\begin{equation}\label{I2def}
I^{2}_{\mu}(\nu)=\begin{cases}
H(\nu \mid \nu_1\otimes\mu)\, ,
& \mbox{ if $\nu\in\skrim_1^*(B)$,} \\
\infty \, ,& \mbox{ otherwise.} \\
\end{cases}
\end{equation}
It is known that $I_\mu^{2}$ is lower semicontinuous and convex.

By $G=G_B\colon \R^d\times\R^d\to [0,\infty]$ we denote the Green function
of a Brownian motion stopped when reaching $\partial B$. That is,
\begin{equation}\label{Greendef}
G(x,y)=\int_0^\infty p_s(x,y) \, \d s,\,\mbox{ for }x,y\in B.
\end{equation}
where $p_s(x,y)$ denotes the transition sub-probability density of the 
stopped motion. Define a function $\skrig\colon\skrim_1(B) \longrightarrow\R$  
by \begin{equation}\label{Gdef}
\skrig(\mu)= \inf_{\heap{\nu\in\skrim^*_1(B)}{\bar{\nu}=\mu}}
\bigl\{ I^{2}_\mu(\nu)-\langle\nu,\log G\rangle \bigr\}
=\inf_{\heap{\nu\in\skrim^*_1(B)}{\bar{\nu}=\mu}}
\int_B\int_B\nu(\d x\,\d y)\log\frac{\nu(\d x\,\d y)}
{\mu(d x)\mu(\d y)G(x,y)}\, .
\end{equation}
Observe that it suffices to take the infimum over measures $\nu$
satisfying $\nu \ll \mu\otimes\mu$.
We can replace $I^{2}_\mu(\nu)$ in the definition of $\skrig$ by either
the relative entropy $H(\nu\,|\,\mu\otimes\mu)$ or the \emph{mutual information}
$H(\nu\,\,|\nu_1\otimes\nu_2)$. 

\subsection{Moment asymptotics}\label{subsec-momasy}

Let $\phi=(\phi_1, \ldots, \phi_n)$ be a family of bounded nonnegative 
Borel measurable functions with compact support in $B$ and define, for any 
$\lambda\in{\mathbb S}_n$,
\begin{equation}\label{ochdef}
\och(\phi, \lambda):=-\inf 
\Big\{ \sum_{i=1}^n \lambda_i H\big( \mu_i \, | \, \phi_i^{2p} \big)
+ p \skrig\Big( \sum_{j=1}^n \lambda_j \mu_j\Big) \, \Big| \, 
\mu_1, \ldots, \mu_n\in \skrim_1(B) \Big\}.
\end{equation}
Note that it suffices to take the infimum over measures $\mu_i$
that satisfy $\mu_i\ll\phi_i^{2p} \, dx$. The following proposition 
is the main result of Section~\ref{sec-momasy}.

\begin{prop}[Asymptotics for mixed moments]\label{mixedmomasy}
Fix $x\in B^p$.
\begin{itemize}
\item[(i)] For every $\lambda\in{\mathbb S}_n$,
\begin{equation}\label{mix1}
\liminf_{k \uparrow \infty} \frac 1k\log\E_x\Bigl[\frac{1}{k!^{p}}
\prod_{i=1}^n\langle\phi_i^{2p},\ell \rangle^{k\lambda_i}\Bigr] \ge 
\och\big(\phi,\lambda\big).
\end{equation}
\item[(ii)] If each $\phi_i$ is bounded away from zero on its support, then
\begin{equation}\label{mix2}
\limsup_{k \uparrow \infty} \sup_{\lambda\in{\mathbb S}_n}
\bigg[ \frac 1k\log \E_x\Bigl[\frac{1}{k!^{p}}\prod_{i=1}^n
\langle\phi_i^{2p},\ell\rangle^{k\lambda_i}\Bigr] - \och\big(\phi,\lambda\big)
\bigg] \le 0.
\end{equation}
\item[(iii)] 
Let $U \subset B$ be a bounded open Lebesgue continuity set whose closure 
is contained in $B$. Let $H\subset \skrim_1(U)$ and $\eps>0$. Then
$$\limsup_{k \uparrow \infty} \frac 1k\log\E_x\Bigl[\frac{1}{k!^{p}} \ell(U)^k 
\1_{\{L\in H\}}\Bigr]  \le -\inf_{\mu \in H[\eps]}
\big\{ I(\mu) + p \skrig(\mu) \big\},$$
where $H[\eps]\subset \skrim_1(U)$ denotes the $\eps$-neighbourhood of the set $H$.
\end{itemize}
\end{prop}

In this section we show how to use this result, together with an analysis of the 
variational problems to complete the proofs of our theorems. We first derive
from Proposition~\ref{mixedmomasy} the moment asymptotics of the random variable 
$\sum_{i=1}^n \langle \phi_i^{2p}, \ell \rangle^{1/p}$. For this purpose define
\begin{equation}\label{Wdef}
W(\phi) =-\inf_{\lambda\in{\mathbb S}_n} \Big\{ 
p I(\lambda) - \och( \phi, \lambda) \Big\}.
\end{equation}

\begin{prop}\label{asym}
Fix $x\in B^p$. 
\begin{enumerate}
\item \ \vspace{-0.8cm}
\begin{equation}\label{lowW}
\liminf_{k\uparrow\infty} \, \frac 1k \, \log \me_x\bigg[ \frac{1}{k!^p} \Big( 
\sum_{i=1}^n \langle \phi_i^{2p}, \ell \rangle^{1/p} \Big)^{kp}\bigg] 
\geq W(\phi).
\end{equation}

\item Assume that every $\phi_i$ is bounded away from zero on its support. Then
\begin{equation}\label{uppW}
\limsup_{k\uparrow\infty} \, \frac 1k \, \log \me_x\bigg[ \frac{1}{k!^p} \Big( 
\sum_{i=1}^n \langle \phi_i^{2p}, \ell \rangle^{1/p} \Big)^{kp}\bigg] 
\leq W(\phi).
\end{equation}
\end{enumerate}
\end{prop}

\begin{Proof}{Proof.} We use Proposition~\ref{mixedmomasy} and denote 
$\varphi_i=\phi_i^{2p}$. For $l,n\in\N$, let
$${\mathbb S}_n(l)=\big\{ \lambda\in {\mathbb S}_n \, : \, l\lambda_i
\mbox{ is an integer for all }i\big\}.$$
The multinomial theorem yields that,
\begin{equation}\label{multinom}
\Big( \sum_{i=1}^n \langle \varphi_i, \ell \rangle^{1/p} \Big)^{kp}
= \sum_{\lambda\in\Sn(kp)} \binom{kp}{kp\lambda_1,\dots,kp\lambda_n}
\prod_{i=1}^n \langle \varphi_i, \ell\rangle^{k\lambda_i}.
\end{equation}
Stirling's formula yields that, uniformly in $\lambda\in {\mathbb S}_n(kp)$, 
as $k\uparrow\infty$,
\begin{equation}\label{IStir}
\binom{kp}{kp\lambda_1,\dots,kp\lambda_n}=e^{-kp I(\lambda)}e^{o(k)}.
\end{equation}

To prove part (i) of Proposition~\ref{asym}, pick some small $\eta>0$ and choose 
an approximate minimiser $\lambda^*\in{\mathbb S}_n$ in \eqref{Wdef} such that 
$$p I(\lambda^*) - \och(\phi, \lambda^*)\leq -W(\phi)+\eta.$$ 
Pick $\lambda\in{\mathbb S}_n(kp)$ such that $|\lambda_i-\lambda^*_i|\leq 
\frac 1k$ for any $i\in\{1,\dots,n\}$. Note that the vector $\widetilde \lambda
=\frac 1\eps(\lambda^*-\lambda(1-\eps))$ lies in ${\mathbb S}_n$.  
Fix $\eps>0$ and use H\"older's inequality to derive
\begin{equation}\label{momlow1}
\E_x\Bigl[\prod_{i=1}^n \langle \varphi_i,\ell\rangle^{k\lambda_i}\Bigr]\geq \E_x\Bigl[\prod_{i=1}^n \langle \varphi_i,\ell\rangle^{k\lambda^*_i}\Bigr]^{\frac 1{1-\eps}} \E_x\Bigl[\prod_{i=1}^n \langle \varphi_i,\ell\rangle^{k\widetilde\lambda_i}\Bigr]^{-\frac \eps{1-\eps}}.
\end{equation}
The last term is further estimated using H\"older's inequality by
\begin{equation}\label{momlow2}
\E_x\Bigl[\prod_{i=1}^n \langle \varphi_i,\ell\rangle^{k\widetilde\lambda_i}\Bigr]^{-\frac \eps{1-\eps}}\geq\prod_{i=1}^n\E_x\bigl[\langle \varphi_i,\ell\rangle^k\bigr]^{-\eps\widetilde\lambda_i/(1-\eps)}\geq \bigl(k!^p e^{Ck}\bigr)^{-\frac \eps{1-\eps}},
\end{equation}
where $C$ only depends on $p$ and $\varphi_1,\dots,\varphi_n$. 
This follows from rough a priori estimates based on Le~Gall's moment formula, 
see e.g.~\cite[Proposition 3]{LG87}. By continuity and nonpositivity of $I$, 
we may estimate $I(\lambda)\leq I(\lambda^*)\frac {1-2\eps}{1-\eps}$, for
$k$ sufficiently large. 
We take expectations on both sides of \eqref{multinom}, and obtain a lower bound by 
restricting the sum on the right side to the unique summand for $\lambda$
chosen above, obtaining, with the help of Proposition~\ref{mixedmomasy}(i), 
$$\begin{aligned}
\liminf_{k\uparrow\infty}& \, \frac 1k \, \log \me_x\bigg[ \frac{1}{k!^p} \Big( 
\sum_{i=1}^n \langle \phi_i^{2p}, \ell \rangle^{1/p} \Big)^{kp}\bigg] \\
& \geq -p I(\lambda)+ \frac 1{1-\eps} \liminf_{k\uparrow\infty}\, \frac 1k \, 
\log \me_x\bigg[ 
\frac{1}{k!^p} \prod_{i=1}^n \langle \phi_i^{2p}, \ell \rangle^{k\lambda_i^*}
\bigg] - \frac{\eps C}{1-\eps}\\
& \geq -\frac 1{1-\eps}\bigl[pI(\lambda^*)-\och(\phi,\lambda^*)
+\eps C-2p\eps I(\lambda^*)\bigr] \\
& \geq-\frac 1{1-\eps}\bigl[ -W(\phi)+\eta+\eps C-2p\eps I(\lambda^*)\bigr].
\end{aligned}$$
Now let $\eta\downarrow0$ and $\eps\downarrow0$ to arrive at \eqref{lowW}.

Now we prove Proposition~\ref{asym}(ii). For any small $\delta>0$ we have, 
by Proposition~\ref{mixedmomasy}(ii), for sufficiently large $k$ and 
any $\lambda\in{\mathbb S}_n$,
$$ \frac 1k\log\E_x\Bigl[\frac 1{k!^{p}}\prod_{i=1}^n\langle\varphi_i,
\ell\rangle^{k\lambda_i}\Bigr] \leq \och(\phi,\lambda)+\delta.$$
Hence, by \eqref{multinom}, \eqref{IStir}, and using that the cardinality of 
$\Sn(kp)$ is bounded by $(kp)^n=e^{o(k)}$, we get
\begin{equation}\label{uppboundfertig}
\limsup_{k\uparrow\infty} \frac 1k \log \me_x\bigg[ \frac{1}{k!^p} \Big( 
\sum_{i=1}^n \langle \phi_i^{2p}, \ell \rangle^{1/p} \Big)^{kp}\bigg] 
\leq -\inf_{\lambda\in\Sn}\bigl[pI(\lambda)-\och(\phi,\lambda)\bigr]+\delta.
\end{equation}
Now we let $\delta\downarrow 0$ in 
\eqref{uppboundfertig} to complete the proof. \end{Proof}\qed

At this point we would like to replace the cumbersome expression~\eqref{Wdef}
for $W(\phi)$ with a more elegant expression involving energies of functions.
The following result is the main result of Section~\ref{sec-anaform}.

\begin{prop}\label{propWident}
Let $\Theta(\phi)$ be as in \eqref{Theta}. Then 
$$W(\phi)=-p \log \frac{\Theta(\phi)}{p}.$$
Moreover, in the variational problems \eqref{Theta} and \eqref{Wdef}
(substituting \eqref{ochdef}) minimisers exist, and they are related
by the formulas
\begin{equation}\label{relmin}
\lambda_i=\|\psi\phi_i\|_{2p}^2, \qquad\mbox{ and }\qquad
\psi^{2p} \phi_i^{2p} = \lambda_i^p \,\frac{d\mu_i}{dx},\qquad
\mbox{ for } i=1,\ldots, n.
\end{equation}
\end{prop}

Using this, we arrive at the following theorem. 

\begin{theorem}[Moment asymptotics]\label{momasy}
Let $\phi=(\phi_1, \ldots, \phi_n)$ be a family of bounded nonnegative Borel 
measurable functions with compact support in $B$, and let $\Theta(\phi)$ be 
as in \eqref{Theta}. Then, for all $x\in B^p$, 
\begin{equation}\label{momasyform}
\lim_{k\uparrow\infty}\frac 1k\log \E_x\Bigl[\frac 1{k!^p}\Big(\sum_{i=1}^n 
\langle \phi_i^{2p}, \ell \rangle^{1/p}\Big)^{kp}\Bigr]=-p\log \frac {\Theta(\phi)}p.
\end{equation}
\end{theorem}

\begin{Proof}{Proof.}
By Propositions~\ref{asym} and~\ref{propWident}, we have the result if 
$\phi_1,\dots,\phi_n$ are bounded away from zero on their 
respective supports. In order to remove this restriction, let $\phi_1,\dots,\phi_n$
be as in the statement and define $\phi_i^{\ssup{\eps}}=\phi_i+
\eps\1_{\,\supp(\phi_i)}$ for $\eps>0$. All we have to show is 
$$\limsup_{\eps\downarrow 0} \Theta(\phi^{\ssup{\eps}})\geq \Theta(\phi).$$ 
It is clear that the condition in the variational problem 
\eqref{Theta} may be relaxed to
$$\Theta(\phi) = \inf\Big\{ \frac p2 \|\nabla \psi\|^2_2 
\, : \, \psi\in\skrid(B), \, \sum_{i=1}^n \| \phi_i \psi \|_{2p}^2 \geq 1 \Big\}.$$
Now we argue that there is $M>0$ such that, for any $\eps\in[0,1]$,
\begin{equation}\label{Thetaeps}
\Theta(\phi^{\ssup{\eps}}) =
\inf\Big\{ \frac p2 \|\nabla \psi\|^2_2 \, : \, \psi\in\skrid(B), \,
\sum_{i=1}^n \| \phi^{\ssup{\eps}}_i \psi \|_{2p}^2 \geq 1, \|\1_U \psi\|_{2p}\leq M \Big\},
\end{equation}
where $U$ denotes the union of the supports of $\phi_1,\dots,\phi_n$. 
Indeed, note that $\eps\mapsto \Theta(\phi^{\ssup{\eps}})$ is decreasing 
and therefore bounded on $[0,1]$. Hence, in dimensions $d\geq 3$, 
\eqref{Thetaeps} follows from \eqref{HoelderRd}, in $d\le 2$ it follows from
\eqref{Hoelderd=2}.

For any $\eps\in[0,1]$ and any $\psi$ in the set on the right hand side of 
\eqref{Thetaeps}, we have
$$\| \phi^{\ssup{\eps}}_i \psi \|_{2p}^2 \leq \| \phi_i \psi \|_{2p}^2 +\eps C,$$
where $C>0$ only depends on $M$ and on the suprema of $\phi_1,\dots,\phi_n$. 
Now choose $\eps\in[0,1/(2C)]$, then we have
$$\begin{aligned}
\Theta(\phi^{\ssup{\eps}})&\geq\inf\Big\{ \frac p2 \|\nabla \psi\|^2_2 
\, : \, \psi\in\skrid(B), \, \sum_{i=1}^n \| \phi_i \psi \|_{2p}^2 \geq 1-\eps C, 
\|\1_U \psi\|_{2p}\leq M \Big\}\\
&\geq\inf\Big\{ \frac p2 \|\nabla \psi\|^2_2 \, : \, \psi\in\skrid(B), \,
\sum_{i=1}^n \| \phi_i \psi \|_{2p}^2 \geq 1-\eps C\Big\}= 
\frac{\Theta(\phi)}{1-\eps C},\end{aligned}$$
which completes the proof.
\end{Proof}\qed

\subsection{Completion of the proofs}\label{sec-finish}

With this result at hand we can easily complete the proof of our main results.
The relation between large moment asymptotics and exponential moments is given
by the following easy lemma.

\begin{lemma}\label{Tauber} Fix $p>0$ and let $X$ be a
nonnegative random variable such that
$$\lim_{k\uparrow\infty}\frac 1k\log E\Bigl[\frac 1{k!^p}X^{kp}\Bigr] 
=-p\log \frac \Theta p$$
exists for some $\Theta>0$. Then we have
$E\bigl[e^X\bigr]<\infty$ if $\Theta>1$, and
$E\bigl[e^X\bigr]=\infty$ if $\Theta<1$.
\end{lemma}

\begin{Proof}{Proof.}
The assumption and Stirling's formula imply that, as $k\uparrow\infty$,
\begin{equation}\label{Stirlappr}
E\bigl[X^{kp}\bigr]=k!^p \Bigl(\frac p\Theta \Bigr)^{kp} e^{o(k)}=(kp)!\Theta^{-pk} e^{o(k)}.
\end{equation}
In the case $\Theta<1$, this shows that the terms $E[X^m]/m!$ 
are not summable over the subsequence $m=kp$, $k\in\N$, 
which implies the second statement. In the case $\Theta>1$ we estimate, 
for any $m$ such that $kp\leq m\leq (k+1)p$, with the help of H\"older's inequality 
and Stirling's formula,
$$E[X^m]\leq E[X^{kp}]^{m/(kp)}=k!^{m/k}\Bigl(\frac p\Theta 
\Bigr)^m e^{o(m)}\leq \Bigl(\frac m{pe}\Bigr)^m \Bigl(
\frac p\Theta \Bigr)^m e^{o(m)}=m!\Theta^{-m}e^{o(m)}.$$
The first statement follows by summing over all $n$.
\end{Proof}\qed

\begin{Proof}{Proof of Theorem~\ref{moments}.}
Apply Lemma~\ref{Tauber} to the situation of Theorem~\ref{momasy}
to get \eqref{criterion}. Similarly to Lemma~\ref{Tauber}, 
\cite[Lemma 2.3]{KM02} relates large integer moments and upper tail 
asymptotics, and allows to infer \eqref{tailbehaviour} from 
Theorem~\ref{momasy}.
\end{Proof}\qed

\begin{Proof}{Proof of Theorem~\ref{LLM}.}
Suppose $H\subset\skrim_1(U)$ is such that $H[\eps]$, an $\eps$-neighbourhood 
of $H$, is disjoint from $\mathfrak{M}$. By~\eqref{relmin}, in the special
case $n=1$, $\phi_1=\1_U$, the set $\mathfrak{M}$ is equal to the set of minimisers 
of $I+p\skrig$ over the set $\skrim_1(U)$. 
As the set of minimisers is closed, we can find $\eps>0$ 
such that $I(\mu)+p\skrig(\mu)-\inf\{I+p\skrig\}>\eps$ for all $\mu\in H[\eps]$. 
Now,
$$\frac 1{a^{1/p}}\log \prob\{ L\in H\, | \,  \ell(U)>a \} = 
\frac 1{a^{1/p}} \log \prob\{  \ell(U)>a \, | \, L\in H\} - 
\frac 1{a^{1/p}} \log \prob\{ \ell(U)>a \},$$
The second term converges to $-p \exp \big( \frac 1p \inf\{I+p\skrig\}\big)$ 
by [Theorem 1.1 and Proposition 2.1, KM02].
For the first term we use the Tauberian Theorem \cite[Lemma 2.3]{KM02},
together with Proposition~\ref{mixedmomasy}(iii) to obtain
\begin{align*}
\lim_{a\uparrow\infty}& \frac 1{a^{1/p}} \log \prob\{ \ell(U)>a \, | \, L\in H\}
\le - p \exp\Big( - \sfrac 1p \lim_{k\uparrow\infty} \frac 1k \log \me_x\Big[ 
\frac{1}{k!^p}\,\ell(U)^k\, \Big| \, L\in H \Big] \Big)\\
& = - p \exp\Big( - \sfrac 1p \lim_{k\uparrow\infty} \frac 1k \log \me_x\Big[ 
\frac{1}{k!^p}\,\ell(U)^k \1_{\{L\in H\}} \Big] \Big) 
\le - p \exp \Big(\sfrac 1p \inf_{\mu\in H[\eps]} \big\{I(\mu)+p\skrig(\mu) 
\big\} \Big).
\end{align*} Altogether, \begin{align*}
\limsup_{a\uparrow\infty}& \frac 1{a^{1/p}} 
\log \prob\{ L\in H\, | \,  \ell(U)>a \}\\ & \le 
 -p \exp\Big(\sfrac 1p\inf_{\mu\in H[\eps]} 
\big\{ I(\mu)+p\skrig(\mu)\big\} \Big)
+p \exp\Big( \sfrac 1p\inf\{I+p\skrig\}\Big)<0,
\end{align*}
which implies the result.\end{Proof}\qed

\section{Moment asymptotics}\label{sec-momasy}

In this section we prove Proposition~\ref{mixedmomasy}. The proof is an 
extension of the proof of \cite[Theorem~1.1]{KM02}, and we are able to 
use some material from there. To keep the notation manageable we assume
that $x=(0,\ldots,0)$, i.e., all motions are started in the origin. We
write $\E$ instead of $\E_0$. The case of arbitrary starting points does 
not pose any additional difficulties for our asymptotic statements. 
We define $\varphi_i=\phi_i^{2p}$. 

The proof is based on a moment formula of Le Gall~\cite{LG86}. To
formulate the result in the necessary generality, recall the Green function
$G=G_B$ and define the function $\Phi_k\colon B^k\to\R$ by
\begin{equation}\label{Phidef}
\Phi_k(y)=\frac 1 {k!} \sum_{\sigma\in\Sym_k}\prod_{i=1}^k
G(y_{\sigma(i-1)},y_{\sigma(i)}),\mbox{ for } y=(y_1,\dots,y_k)\in\ B^k,
\end{equation}
where we put $y_0=0$, the starting point of the motions.
$\Sym_k$ is the symmetric group in $k$ elements and we write
elements $\sigma\in\Sym_k$ as permutations $\sigma\colon\{1,\dots,k\}\to
\{1,\dots,k\}$ and agree on the additional convention that $\sigma(0)=0$.
Introduce the \emph{empirical measure} of the vector $y=(y_1,\dots,y_k)\in
B^k$, $$L_{y,k}=\frac 1{k}\sum_{j=1}^{k} 
\delta_{y_j}\in\skrim_1(B),$$
and note that $\Phi_k$ is a permutation symmetric function, i.e., it depends
on $y$ only via $L_{y,k}$. $\Phi_k$ assumes the value $\infty$ if and only 
if $0,y_1,\dots,y_k$ are not pairwise distinct. For the rest of the proof,
we tacitly assume that the numbers $k\lambda_1,\dots, k\lambda_n$ are integers.
This simplification can be justified by simple local approximations of the
type $\lambda_i \leadsto \lfloor k\lambda_i \rfloor/k$ or
$\lambda_i \leadsto \lceil k\lambda_i \rceil/k$. Inequalities~\eqref{momlow1}
and~\eqref{momlow2} show that this does not affect the asymptotics on the left
hand side of \eqref{mix1} and \eqref{mix2}.

In the following, we organise the vector $y$ as 
$$y=(y^i_j \, : \, j=1,\dots, k\lambda_i, \mbox{ and } i=1,\ldots,n).$$ 
\begin{lemma}[Moment formulas]\ \\ 
\vspace{-0.7cm}\\
\begin{enumerate}
\item \ \vspace{-0.9cm}
\begin{equation}\label{momentformula}
\me\bigg[ \frac{1}{k!^p} \prod_{i=1}^n \langle \varphi_i, \ell \rangle^{k\lambda_i} 
\bigg]=\int_{B^k}dy\,\big(\Phi_k(y)\big)^p\, \prod_{i=1}^{n}
\prod_{j=1}^{\lambda_ik} \varphi_i(y^i_j).
\end{equation}
\item Let $H\subset\skrim_1(U)$ be a Borel set, and $\delta>0$. 
Then, for all sufficiently large $k\in\N$, 
\begin{equation}\label{genmomentformula}
\me\Big[ \frac{1}{k!^p}\, \ell(U)^k \1_{\{ L \in H\}}\Big] \le
2 \int_{U^k} dy\,\1_{\{L_{y,k}\in H[\delta]\}}\,
\big(\Phi_k(y)\big)^p,\end{equation}
where $H[\delta]$ is the $\delta$-neighbourhood of $H$.
\end{enumerate}
\end{lemma}

\begin{Proof}{Proof.} A variant of part~(i) was proved in \cite{LG86}, 
see~\cite[(2c)]{LG87}. We focus on part~(ii), where the main idea is that 
the random variable $L$  is asymptotically estimated in terms of the 
empirical measure $L_{y,k}$ induced by the integration variable 
$y=(y_1,\dots,y_k)$ on the right hand side of the moment formula. 

An essential ingredient is Le Gall's Wiener sausage characterization of 
$\ell$, see \cite[Th.~3.1]{LG86}. For every $\eps>0$ define the 
{\it Wiener sausage\/} around $W_{i}$ by 
\begin{equation}\label{sausage}
S_\eps^i=\Big\{ x\in B \, : \, \mbox{ there is } t\in[0,T_i) \mbox{ with }
|x-W_{i}(t)|<\eps \Big\},\mbox{ for }i=1,\dots,p \, ,
\end{equation}
and their intersection $S_\eps=\bigcap_{i=1}^p S_\eps^i$.
Recall that $S=W^1[0,T^1]\cap\ldots\cap W^p[0,T^p]$
and observe that $S=\bigcap_{\eps>0} S_\eps$ 
is the intersection of the $p$ independent Brownian paths. 
Define
\begin{equation}\label{sdef}
s_d(\eps)=
\begin{cases}
\pi^{-p}\log^p(1/\eps), &
\mbox{if $d=2$,}\\
(2\pi\eps)^{-2}, & 
\mbox{if $d=3$ and $p=2$,}\\
\frac{2}{\omega_d (d-2)}\eps^{2-d}, & 
\mbox{if $d\geq 3$ and $p=1$,}\\
\end{cases}
\end{equation}
where $\omega_d$ is the volume of the $d$-dimensional unit ball. 
For every set $A\subset B$ that is almost surely an $\ell$-continuity set,
\begin{equation}\label{Legallformula}
\lim_{\eps\downarrow 0} s_d(\eps) \Leb(S_\eps \cap A) = \ell(A)\, ,
\end{equation}
in the $L^k(\prob)$-sense for all positive integers $k$ and, in particular, 
in probability. We denote by 
$\Lambda_{\eps}$ the normalized restriction of the 
Lebesgue measure to $S_\eps\cap U$, considered as a probability measure
on $U$, {i.e.,} as an element of $\skrim_1(U)$. As we have assumed that
$U$ is a Lebesgue continuity set, we see from the moment 
formula~\eqref{momentformula} 
for $n=k=1$ and $\varphi_1=\1_{\partial U}$ that almost surely 
$\ell(\partial U)=0$ and thus (\ref{Legallformula}) applies to the set $A=U$. 
This implies in particular 
that for the uniform distributions on the intersections of the Wiener sausage,
\begin{equation}\label{probconv}
\lim_{\eps \downarrow 0} \Lambda_{\eps} = L \mbox{ in probability.}
\end{equation}
Clearly, for every $H\subset\skrim_1(U)$ and every positive integer $k$, 
the random variables $s_d(\eps) \Leb  (S_\eps\cap U)\1_{\{L \in H\}}$ 
converge as $\eps\downarrow 0$ in the $L^k(\prob)$-sense to 
$\ell(U)\1_{\{L \in H\}}$. Thus
\begin{equation}\label{firststep}
\me\Big[ \ell(U)^k \1_{\{L \in H\}}\Big]  =
\lim_{\eps\downarrow 0} s_d(\eps)^{k} 
\me\Big[\Leb( S_\eps\cap U)^k \1_{\{L\in H\}}\Big]\, .
\end{equation}

Suppose now that $\delta>0$ and a Borel set $H\subset\skrim_1(\R^d)$ are given. 
Using (\ref{firststep}) we obtain for the left hand side in our statement,
\begin{equation}
\begin{aligned}
\me\Big[ \ell(U)^k \1_{\{L \in H\}}\Big]  = &
\lim_{\eps\downarrow 0} s_d(\eps)^{k} \,
\me\Big[\Leb( S_\eps\cap U)^k \1_{\{L\in H\}}\Big] \\
\leq  & \limsup_{\eps\downarrow 0} s_d(\eps)^{k} \,
\me\Big[\int_{(U\cap S_\eps)^k} dy \, 
\1_{\{L_{y,k}\in H[\delta]\}}\Big] \label{sum1}\\
& + \limsup_{\eps\downarrow 0} s_d(\eps)^{k} \,
\me\bigg[\Leb( S_\eps\cap U)^k {\tt E}_{\Lambda_{\eps}^{\otimes k}}\Big[ 
\1_{\{L\in H\}}-\1_{\{L_{Y,k}\in H[\delta]\}}\Big]\bigg] 
\, .
\end{aligned}
\end{equation}
Here ${\tt E}_{\Lambda_{\eps}^{\otimes k}}$ denotes expectation with respect to 
$\Lambda_{\eps}^{\otimes k}$, and $Y\colon (S_\eps\cap U)^k\to (S_\eps\cap U)^k$ is 
the identity map, i.e., a $(S_\eps\cap U)$-valued random variable
with distribution $\Lambda_{\eps}^{\otimes k}$.

To treat the first term on the right of \eqref{sum1} we recall from \cite{LG86} 
that the family of functions
$$y \mapsto s_d(\eps)^{k} \, \prob\Big\{
y_j\in S_\eps\mbox{ for all } j \Big\},\qquad 
\mbox{ for } \eps\in(0,{\textstyle{\frac 12}}),$$
is dominated by an integrable function and,  
$$\lim_{\eps\downarrow 0} s_d(\eps)^{k} \, \prob\Big\{
y_j\in S_\eps\mbox{ for all } j \Big\}
=\Big( \sum_{\sigma\in \Sym_k} \prod_{i=1}^k G\Big(y_{\sigma(i-1)}, y_{\sigma(i)} 
\Big) \Big)^p\, .$$
Hence for the first term on the right hand side of (\ref{sum1}), 
by Fubini's Theorem and the theorem of dominated convergence,
\begin{equation}
\begin{aligned}
\lim_{\eps\downarrow 0} \, s_d(\eps)^{k} \,
\me\Big[\int_{(U\cap S_\eps)^k}dy\, 
\1_{\{L_{y,k}\in H[\delta]\}}\Big]
& = \lim_{\eps\downarrow 0} \int_{U^k}dy\, 
\1_{\{L_{y,k}\in H[\delta]\}} s_d(\eps)^{k} \, \prob\Big\{
y_j\in S_\eps\mbox{ for all } j \Big\}\\
& = \int_{U^k}  dy\, 
\1_{\{L_{y,k}\in H[\delta]\}}  \lim_{\eps\downarrow 0} s_d(\eps)^{k} 
\, \prob\Big\{ y_j\in S_\eps\mbox{ for all } j \Big\} \\
& = \int_{U^k}dy\, 
\1_{\{L_{y,k}\in H[\delta]\}} \Big( \sum_{\sigma\in \Sym_k} 
\prod_{i=1}^k G\Big(y_{\sigma(i-1)}, y_{\sigma(i)} \Big) \Big)^p \, ,
\label{formula-a}\end{aligned}\end{equation}
which is equal to the right hand side in our claim \eqref{genmomentformula}.

It remains to derive an upper bound  for the second term of the right hand side of
\eqref{sum1} that is negligible with respect to the left hand side in 
\eqref{genmomentformula}. Observe that, 
\begin{equation}
\begin{aligned}
\1_{\{L\in H\}} -\1_{\{L_{Y,k}\in H[\delta]\}}
&\leq \1_{\{L\in H\}}\1_{\{L_{Y,k}\notin H[\delta]\}}\\
&\leq \1_{\{d(L,\Lambda_{\eps})\geq \delta/2\}}+\1_{\{L\in H\}}
\1_{\{L_{Y,k}\notin H[\delta]\}}\1_{\{d(L,\Lambda_{\eps})< \delta/2\}}\\
&\leq \1_{\{d(L,\Lambda_{\eps})\geq \delta/2\}}+\1_{\{L\in H\}}
\1_{\{d(\Lambda_{\eps},L_{Y,k})\geq \delta/2\}}.
\end{aligned}
\end{equation}
Hence, we obtain,
\begin{equation}\begin{aligned}
\limsup_{\eps\downarrow 0}  s_d(\eps)^{k} & \,
\me\Big[\Leb( S_\eps\cap U)^k {\tt E}_{\Lambda_{\eps}^{\otimes k}}\Big[ 
\1_{\{L\in H\}}-\1_{\{L_{Y,k}\in H[\delta]\}}\Big]\Big] \\ 
\leq & \limsup_{\eps\downarrow 0} s_d(\eps)^{k} \,
\me\Big[\Leb( S_\eps\cap U)^k 
\1_{\{d(L,\Lambda_{\eps})\geq \delta/2\}}\Big] \\
&  + \limsup_{\eps\downarrow 0} s_d(\eps)^{k} \,
\me\Big[\Leb( S_\eps\cap U)^k \1_{\{L\in H\}} \Lambda_{\eps}^{\otimes k}\Big\{
d(\Lambda_{\eps},L_{Y,k})\geq \delta/2 \Big\} \Big] \label{sum4}\, .
\end{aligned}
\end{equation}
Recall that the family $s_d(\eps)\Leb (S_\eps\cap U)$ with
$0<\eps<1/2$ is bounded in every $L^k(\prob)$. Using this fact together
with (\ref{probconv}), we get
\begin{equation}\label{formula-b}
\limsup_{\eps\downarrow 0} s_d(\eps)^{k} \,
\me\Big[\Leb( S_\eps\cap U)^k \1_{\{d(L,\Lambda_{\eps})\geq \delta/2\}}\Big] = 0 \, .
\end{equation}
To show that the remaining summand on the right hand side of (\ref{sum4}) is small, 
we need an upper bound for 
$\Lambda_{\eps}^{\otimes k}\{d(\Lambda_{\eps},L_{Y,k})\geq \delta/2\}$,
which depends neither on $\eps$ nor on the Brownian paths. To achieve
this, we first use the fact that the weak topology of measures can be 
approximated by a finite partition. More precisely, for the given $\delta>0$, 
we can find a finite partition $\skrip$ of $B$ and an $\eta=\eta(\delta)>0$ 
such that,
$$
\sup_{M \in \skrip} |\nu(M)-\mu(M)| \leq \eta\qquad
\text{ implies } \qquad d(\nu,\mu) <\delta/2\,  .
$$
Denote the cardinality of $\skrip$ by $N$ and the canonical projection
by $\pi\colon B\to\{1,\ldots,N\}$ and recall that for any two probability measures
$P$, $Q$ on $\{1,\ldots,N\}$ the Kullback-Leibler distance $H(P\mid Q)$ 
is bounded from below by $\frac 12 \sup_{A\subset\{1,\ldots,N\}}|P(A)-Q(A)|^2$. 
Hence, we obtain for any two probability measures
$\nu$, $\mu$ on $\R^d$ that
$$
d(\nu,\mu)\geq\delta/2 \qquad\text{ implies } \qquad H(\nu\circ\pi^{-1}\mid \mu\circ\pi^{-1})
\geq {\eta^2/2}\, .
$$
Denote by $\Gamma$ the set of probability measures $Q$ on $\{1,\ldots,N\}$
with $H(Q\mid \Lambda_{\eps}\circ\pi^{-1})\geq {\eta^2/2}$. {F}rom the
upper bound in Sanov's Theorem for the finite alphabet $\{1,\ldots,N\}$,
see \cite[p.~15]{DZ98}, we infer that
\begin{equation}\begin{aligned}
\Lambda_{\eps}^{\otimes k}\{ d(\Lambda_{\eps},L_{Y,k})\geq \delta/2\}  
& \leq (k+1)^N \exp\Big(-k \inf_{Q\in\Gamma} H(Q\mid \Lambda_{\eps}
\circ\pi^{-1})\Big) \\
& \leq (k+1)^N \exp(-k{\eta^2/2})\, ,\label{Sanov}
\end{aligned}\end{equation}
which is the required upper bound. We infer from this and (\ref{firststep}) that,
\begin{equation}\begin{aligned}
\limsup_{\eps\downarrow 0} s_d(\eps)^{k} & \,\me\Big[\Leb( S_\eps\cap U)^k 
\1_{\{L\in H\}}\Lambda_{\eps}^{\otimes k}\Big\{ d(L_{Y,k}, \Lambda_{\eps}) 
\geq \delta/2 \Big\} \Big] \\
& \leq \limsup_{\eps\downarrow 0} s_d(\eps)^{k} \,
\me\Big[\Leb( S_\eps\cap U)^k \1_{\{L\in H\}}(k+1)^N \exp(-k{\eta^2/2})\Big]\\
& \leq (k+1)^Ne^{-k{\eta^2/2}}\me\Big[ \ell(U)^k \1_{\{L \in H\}}\Big]\, .
\label{formula-c}
\end{aligned}\end{equation}
Putting (\ref{formula-a}), (\ref{formula-b}) and
(\ref{formula-c}) together, we obtain
$$(1-(k+1)^Ne^{-k{\eta^2/2}})\me\Big[ \ell(U)^k \1_{\{L \in H\}}\Big]
\leq \int_{U^k}dy\, \1_{\{L_{y,k}\in H[\delta]\}}
\Big( \sum_{\sigma\in \Sym_k} \prod_{i=1}^k G\Big(y_{\sigma(i-1)},y_{\sigma(i)} 
\Big) \Big)^p.$$
Finally, $1-(k+1)^N\exp(-k{\eta^2/2})>1/2$ for all sufficiently large $k$, and 
this finishes the proof.\end{Proof}\qed

In order to derive the upper bounds in Proposition~\ref{mixedmomasy}
it is  necessary to replace the Green function $G$ in the definition 
of $\Phi_k(y)$ by some bounded function. We achieve this by cutting off 
at a large level and show that this does not change the exponential 
rate of $\Phi_k(y)$ asymptotically as the cut-off level gets large. 
Introduce, for $M\geq 0$, the cut-off Green function 
$G_M=G\wedge M$ and denote,
\begin{equation}\label{GMexpression}
\Phi_{k,M}(y)=\frac 1{k!}\sum_{\sigma\in \Sym_k} 
\prod_{i=1}^k 
G_M\left(y_{\sigma(i-1)}, y_{\sigma(i)} \right), \quad \mbox{ for } 
y\in B^k.
\end{equation}

The following lemma provides the cutting argument:

\begin{lemma}\label{cutoff}
There is $C_0>0$ and, for all sufficiently large $M>1$ and small $\eta\in(0,1)$, 
there are constants $C_M>0$ and $\eps_\eta>0$ such that 
$\lim_{M\uparrow\infty}C_M=\lim_{\eta\downarrow 0}\eps_\eta=0$,
and the following holds.
\begin{enumerate}
\item For any $k\in\N$ and for any $\lambda\in{\mathbb S}_n$, 
\begin{equation}\label{cutG}
\begin{aligned}
\int_{B^{k}}dy\, & \bigl(\Phi_{k}(y)\bigr)^p\prod_{i=1}^{n}
\prod_{j=1}^{\lambda_ik} \varphi_i(y^i_j)\\
&\leq 2^ppk (2C_0)^k C_M^{\eta k}+ 2^p (1+\eps_\eta)^k p
\sum_{m=\lceil k(1-p\eta)\rceil}^k \int_{B^{m}}dy\, \big(\Phi_{m,M}(y)\big)^p\prod_{i=1}^{n}
\prod_{j=1}^{\lambda_im} \varphi_i(y^i_j).
\end{aligned}
\end{equation}
\item For any $H\subset \skrim_1(U)$ Borel, and $\delta>0$, 
one can pick $\eta>0$ small enough such that 
\begin{equation}\label{cutG2}\begin{aligned}
\int_{U^{k}}dy\, & \bigl(\Phi_{k}(y)\bigr)^p \1_{\{L_{y,k}\in H\}} \\
& \leq 2^ppk (2C_0)^k C_M^{\eta k}+ 2^p (1+\eps_\eta)^k p
\sum_{m=\lceil k(1-p\eta)\rceil}^k \int_{U^{m}}dy\, \big(\Phi_{m,M}(y)\big)^p
\1_{\{L_{y,m}\in H[\delta]\}}.
\end{aligned}
\end{equation}
\end{enumerate}
\end{lemma}

\begin{Proof}{Proof.}
This follows from an obvious adaptation of Lemmas~3.2 and 3.3 in
\cite{KM02} and their proofs (recall that $\varphi_1,\dots,\varphi_n$ are bounded).
\end{Proof}
\qed

Note that for the upper bounds we may now focus on the $m$-fold integral on the 
right hand side of \eqref{cutG}, resp.~\eqref{cutG2}. 
Observe that the integration domain $B^m$ may 
be replaced by the compact set $U^m$ where $U$ is the union of the supports
of $\varphi_i$ for $i=1,\ldots, n$.

Our second main technical tool is a reduction to a discrete counting 
argument. For this purpose, we introduce a finite partition of $U$
which is carefully chosen in order to represent 
many details of the continuous picture. 

To introduce appropriate notation, let $\Sigma_r=\{1,\dots,r\}$ and denote the 
partition sets by $U_1, \dots, U_r$. We assume that every $U_l$ is measurable 
and has positive Lebesgue measure. In Lemma~\ref{fineness} below we shall
make precise how fine we choose the partition. We call 
$\pi\colon U\to \Sigma_r$ the canonical projection, that is, $x\in U_{\pi x}$ 
for any $x\in U$. We write $\pi y=((\pi y^i_1,
\dots,\pi y^i_{\lambda_i m})\, : \, i=1, \dots,n)$
if $y=((y^i_1,\dots,y^i_{\lambda_i m})\, : \, i=1, \dots,n)$. 
If $\mu$ is a probability measure on $U$, then $\pi\mu\in\skrim_1(\Sigma_r)$ 
is its projection on $\Sigma_r$. Similarly for $\nu\in\skrim_1(U^2)$ we 
denote the projection on $\Sigma_r^2$ by $\pi\nu\in\skrim_1(\Sigma_r^2)$.
If $v$ is in the set $\skrim_1^*(\Sigma_r)$ of probability measures on 
$\Sigma_r^2$ with equal marginals, we denote by $\overline v\in\skrim_1(\Sigma_r)$ 
its left or right marginal measure.
Note that $\pi \overline \nu=\overline{\pi\nu}$ for any $\nu\in \skrim_1^*(U)$,
where $\overline\nu$ is the marginal measure of $\nu$.

For measures $u\in\skrim_1(\Sigma_r)$ and
$v\in\skrim_1(\Sigma_r^2)$ we define discrete analogues of the relative 
entropy functionals $I$ and $I^{ 2 }_\mu$ by 
\begin{equation}\label{Itildedef}
\widetilde I(u)= \sum_{l\in \Sigma_r} u_l \, 
\log\frac{u_l}{|U_l|}
\qquad\mbox{ and }\qquad \widetilde I_u^{2}(v)=\sum_{l,m\in \Sigma_r} 
v_{l,m} \, \log\frac{v_{l,m}}{\overline v_l u_m},
\end{equation}
using the usual convention $0\log 0=0$.
Recall that $G_M=G\wedge M$ and define the approximate Green functions 
$G_M^+, G^-\colon \Sigma_r^2\to\R$ by
\begin{equation}\label{Gapprdef}
G_M^+(l,m)=\sup_{\heap{x\in U_l}{y \in U_m}} G_M(x,y)\qquad \mbox{and}
\qquad
G^-(l,m)=\inf_{\heap{x\in U_l}{y \in U_m}} G(x,y).
\end{equation}
Functions $\skrig_{M}^+$ and $\skrig^-$ on $\skrim_1(\Sigma_r)$ analogous 
to $\skrig$ in (\ref{Gdef}) are defined by
\begin{equation}\label{skrigapprdef}
\skrig_{M}^+(u)= \inf_{\heap{v\in\skrim^*_1(\Sigma_r)}{\overline v=u}}
\left\{ \widetilde I^{2}_u(v)-\left\langle v,\log G_{M}^+
\right\rangle\right\}\quad \mbox{ and }\quad
\skrig^-(u)= \inf_{\heap{v\in\skrim^*_1(\Sigma_r)}{\overline v=u}}
\left\{ \widetilde I^{2}_u(v)-\left\langle v,\log G^-
\right\rangle\right\},
\end{equation}
where we used the notation $\langle v,F\rangle = \sum_{l,m\in\Sigma_r} v_{l,m}
\,F(l,m)$. 

The functions $\skrig_{M}^+$ and $\skrig^-$ are continuous.
Indeed, for fixed $u$, if the set $\widetilde V\subset\skrim_1^*(\Sigma_r)$ is 
a neighbourhood of the set $\{v\in\skrim^*_1(\Sigma_r)\, :\, \overline v=u\}$, 
there exists a neighbourhood
$\widetilde U$ of $u$ with $\{v\in\skrim^*_1(\Sigma_r)\, :\, 
\overline v=\tilde{u}\}\subset \widetilde V$
for all $\tilde{u}\in \widetilde U$. Together with the obvious continuity of
$\widetilde I^{2}_u(v)$ in both arguments $u$ and $v$ and of 
$v\mapsto \langle v,F\rangle$ this implies continuity of $\skrig_{M}^+$ and $\skrig^-$.

Introduce the empirical measure of the vector $\sigma^i=(\sigma^i_1,\dots,
\sigma^i_{\lambda_im})\in \Sigma_r^{\lambda_im}$ of length $\lambda_im$ by 
putting
\begin{equation}
L_{\sigma^i,\lambda_im}=\frac 1{\lambda_im}\sum_{j=1}^{\lambda_im} 
\delta_{\sigma^i_j}\in\skrim_1(\Sigma_r),
\end{equation}
and the global empirical measure of $\sigma=(\sigma^1,\ldots, \sigma^n)\in
\Sigma_r^m$ by
\begin{equation}
L_{\sigma,m}=\frac 1{m}\sum_{i=1}^{n}\sum_{j=1}^{\lambda_im} \delta_{\sigma^i_j}
= \sum_{i=1}^n \lambda_i L_{\sigma^i,\lambda_im}\in\skrim_1(\Sigma_r).
\end{equation}

By \cite[Lemma~3.5]{KM02}, for any $M>0$, uniformly in $y\in U^m$, as $m\uparrow \infty$,
\begin{eqnarray}
\Phi_{m,M}(y)&\leq& e^{o(m)}\exp\Bigl(-m\skrig_M^+\bigl(L_{\pi y,m}\bigr)\Bigr),
\label{Phidisc}\\
\Phi_{m}(y)&\geq&e^{o(m)}\exp\Bigl(-m\skrig^-\bigl(L_{\pi y,m}\bigr)\Bigr).
\end{eqnarray}
Now we go back to the integral on the right hand side of \eqref{cutG} and 
rewrite the integral over $B^m$ (which we have replaced by the integral over $U^m$)
as an integral over the partition sets $U{\mbox{\scriptsize $\sigma^i_j$}}$ 
with $i\in\{1,\dots,n\}$ and $j\in\{1,\dots,\lambda_i m\}$ and sum over all 
$\sigma^i_j\in\Sigma_r$.
 
Note that, for any $\sigma\in \Sigma_r^m$, the map $y\mapsto L_{\pi y,m}=
L_{\sigma,m}$ is constant on the set of $y\in U^m$ satisfying 
$y^i_j\in U_{\sigma^i_j}$, where we again organise $\sigma\in \Sigma_r^m$ as 
$\sigma=(\sigma^i_j\,:\,j=1,\dots,\lambda_i m,\, i=1,\dots,n)$. Hence, 
$$
\begin{aligned}
\int_{B^{m}}dy\, \big(\Phi_{m,M}(y)\big)^p\prod_{i=1}^{n}
\prod_{j=1}^{\lambda_im} \varphi_i(y^i_j)
&\leq \sum_{\sigma\in \Sigma_r^m}e^{-pm \skrig_M^+(L_{\sigma,m})}\prod_{i=1}^{n}
\prod_{j=1}^{\lambda_im} \int_{U_{\sigma^i_j}}\varphi_i(y^i_j)\,d y^i_j\\
&=\sum_{\sigma\in \Sigma_r^m}e^{-pm \skrig_M^+(L_{\sigma,m})}\prod_{i=1}^{n}
\exp\Bigl(m\lambda_i \bigl\langle L_{\sigma^i,m},\log \int_{U_\cdot}\varphi_i\bigr\rangle\Bigr).
\end{aligned}
$$
Analogously, we have a lower bound for the left hand side with $\Phi_{m,M}$ replaced
by $\Phi_m$ in terms of the right hand side with $\skrig_M^+$ replaced by $\skrig^-$.
We rewrite the sum over $\sigma\in\Sigma_r^m$ as $n$ sums over probability measures
$u_i\in\skrim_1(\Sigma_r)$ and count the numbers of 
$\sigma^i\in\Sigma_r^{\lambda_i m}$ such that $u_i$ is the empirical 
measure of~$\sigma^i$,
$$
\begin{aligned}
\int_{B^{m}}dy\, &\big(\Phi_{m,M}(y)\big)^p\prod_{i=1}^{n}
\prod_{j=1}^{\lambda_im} \varphi_i(y^i_j)\\
&\leq \sum_{\heap{u_i\in\skrim_1^{(m\lambda_i)}(\Sigma_r)}
{\forall i=1,\dots,n}}e^{-mp\skrig_M^+(\sum_{i=1}^n \lambda_i u_i)}
\Bigl[\prod_{i=1}^n\#\{\sigma^i\in\Sigma_r^{\lambda_i m}\colon u_i=L_{\sigma^i,m}\}\Big]
\prod_{i=1}^n e^{m\lambda_i \langle u_i,\log \int_{U_\cdot}\varphi_i\rangle},
\end{aligned}
$$
where  $\skrim_1^{(m\lambda_i)}(\Sigma_r)$ is the set of those
 $u_i$ such that $m\lambda_i u_i(l)$ is an integer for any $l\in\Sigma_r$.
By Stirling's formula, the $i$th counting factor on the right is equal to
$e^{o(m)}\exp(-m\lambda_i \sum_{l\in\Sigma_r}u_i(l)\log u_i(l))$, uniformly in $u_i\in 
\skrim_1^{(m\lambda_i)}(\Sigma_r)$ and in $\lambda\in{\mathbb S}_n$. Indeed, choose
$0<c<C$ such that $c\leq n!/[(\frac ne)^{-n}n^{-1/2}]\leq C$ for any $n\in\N$, 
and estimate, for $\lambda_i>0$,
$$
\begin{aligned}
\#\{\sigma^i\in\Sigma_r^{\lambda_i m}&\colon u_i=L_{\sigma^i,m}\}
=\frac {(\lambda_im)!}{\prod_{l\colon u_i(l)>0}(u_i(l)\lambda_i m)!}\\
&\leq \exp\Bigl\{-m\lambda_i \sum_{l}u_i(l)\log u_i(l)\Bigr\}Cc^{-r}\sqrt{\frac{\lambda_i}
{\prod_{l\colon u_i(l)>0}(u_i(l)\lambda_i)}}.
\end{aligned}
$$
Now use that $u_i(l)\lambda_i \geq \frac 1m$ if $u_i(l)>0$. The lower bound
is derived in a similar way.

Note that the cardinality of $\skrim_1^{(m\lambda_i)}(\Sigma_r)$ is polynomial 
in $m$, uniformly in $\lambda\in{\mathbb S}_n$. Hence, we obtain
\begin{equation}\label{limplus}
\limsup_{m\uparrow \infty} \sup_{\lambda\in{\mathbb S}_n}
\bigg[ \frac 1m\log\Bigl(\int_{B^{m}}dy\,
 \big(\Phi_{m,M}(y)\big)^p\prod_{i=1}^{n}
\prod_{j=1}^{\lambda_im} \varphi_i(y^i_j)\Bigr)
-\widetilde\och_M^+(\phi,\lambda) \bigg] \leq 0,
\end{equation}
where we recall that $\varphi_i=\phi_i^{2p}$ and introduce
$$
\widetilde\och_M^+(\phi,\lambda)=
-\inf_{u_1,\dots,u_n\in\skrim_1(\Sigma_r)}\Bigl[p\skrig_M^+
\Bigl(\sum_{i=1}^n\lambda_i u_i\Bigr)
+\sum_{i=1}^n \lambda_i H\Bigl(u\,\Big|\,{\int_{U_\cdot}\varphi_i}\Bigr)\Bigr].
$$
It is easy to see that the map $\widetilde\och_M^+(\phi,\cdot)$ is continuous 
on the simplex ${\mathbb S}_n$. Indeed, the family of mappings
$$
\lambda\mapsto p\skrig^+_M\Bigl(\sum_{i=1}^n\lambda_i u_i\Bigr)+\sum_{i=1}^n 
\lambda_i H\Bigl(u_i\Big|\int_{U_\cdot}\varphi_i\Bigr),\qquad 
\mbox{ for } u_1,\dots,u_n\in\skrim_1(\Sigma_r),
$$
is uniformly equicontinuous on ${\mathbb S}_n$, since $\skrig_M^+$ is uniformly 
continuous on $\skrim_1(\Sigma_r)$, and the map $u\mapsto H(u\mid\int_{U_\cdot}\varphi_i)$ 
is bounded and continuous for every $i$.

Analogously to \eqref{limplus}, we have, for any $\lambda\in{\mathbb S}_n$,
\begin{equation}\label{limminus}
\liminf_{k\uparrow \infty}
\frac 1k\log\Bigl(\int_{B^{k}}dy\, \big(\Phi_{k}(y)\big)^p\prod_{i=1}^{n}
\prod_{j=1}^{\lambda_ik} \varphi_i(y^i_j)\Bigr)
\geq \widetilde\och^-(\phi,\lambda),
\end{equation}
where $\widetilde\och^-(\phi,\lambda)$ is defined as 
$\widetilde\och^+_M(\phi,\lambda)$ with $\skrig_M^+$ replaced by
$\skrig^-$. Like $\widetilde\och_M^+(\phi,\cdot)$, the function  
$\widetilde\och^-(\phi,\cdot)$ is also continuous on ${\mathbb S}_n$.

Now we determine the fineness of the partition $(U_1,\dots,U_r)$ of the set $U$. 

\begin{lemma}[Choice of the partition]\label{fineness} 
Fix $M>0$ and any $\delta>0$.
\begin{enumerate}
\item For any $\lambda\in{\mathbb S}_n$, the partition $(U_1,\dots,U_r)$ 
of the set $U$ may be chosen so fine that 
$$\widetilde\och^-(\phi,\lambda)\geq\och(\phi,\lambda)-\delta.$$
\item If each $\phi_i$ is bounded away from zero on its support, then
the partition $(U_1,\dots,U_r)$ of the set $U$
may be chosen so fine that, for any $\lambda\in{\mathbb S}_n$,
$$\widetilde\och^+_M(\phi,\lambda)\leq \och(\phi,\lambda)+\delta.$$
\item For any $H\subset\skrim_1(U)$ and $\delta>0$ the partition
$(U_1,\dots,U_r)$ of the set $U$ may be chosen so fine that 
$$-\inf_{\mu\in H}\big\{\widetilde{I}(\pi\mu) + p \skrig^+_M(\pi\mu) \big\}
\leq - \inf_{\mu\in H[\delta]} \big\{ I(\mu) + p \skrig(\mu) \big\} +\delta.$$
\end{enumerate}
\end{lemma}

\begin{Proof}{Proof.} We first prove (i). Choose approximate minimisers 
$\mu_1^*,\dots,\mu_n^*\in\skrim_1(U)$ in \eqref{ochdef} and an approximate 
minimiser $\nu^*\in\skrim^*_1(B^2)$ in the definition \eqref{Gdef} 
of~$\skrig$ (satisfying $\overline\nu^*=\sum_{i=1}^n\lambda_i\mu_i^*$) such that
\begin{equation}\label{ochapprmin}
-\och(\phi,\lambda)\geq \sum_{i=1}^n \lambda_i H( \mu_i^* \, | \, \varphi_i) 
+p\bigl( I^2_{\overline\nu^*}(\nu^*)-\langle\nu^*,\log G\rangle\bigr)-\frac \delta 2.
\end{equation}
We now use Jensen's inequality to show that, for any partition and for any 
$\mu$ resp.~$\nu$, we have
\begin{equation}\label{HIesti}
H( \mu \mid \varphi_i)\geq H\Bigl(\pi\mu\,\Big|\,\int_{U_\cdot}\varphi_i\Bigr)
\qquad\mbox{and}\qquad I^2_{\overline\nu}(\nu)\geq \widetilde I^2_{\pi\overline\nu}(\pi\nu).
\end{equation}
To prove this, abbreviate $f(y)=y\log y$ and note that
$$H( \mu \mid \varphi_i)=\sum_l\int_{U_l}\varphi_i(x)\,dx\int_{U_l}
\frac{dx\,\varphi_i(x)}{\int_{U_l}\varphi_i(x)\,dx}\,
f\Bigl(\frac{\mu(dx)/dx}{\varphi_i(x)}\Bigr).$$
Now use Jensen's inequality for the convex function $f$ and summarise 
to arrive at the first inequality in \eqref{HIesti}. The other one is 
proved analogously, noting that $I^2_{\overline \nu}(\nu)=\langle 
\overline\nu\otimes\overline\nu, f\circ (d\nu/d(\overline\nu\otimes\overline\nu))
\rangle$.

For $\nu^*$ fixed above, we may choose the partition $(U_1,\dots,U_r)$ 
of $U$ so fine that $\langle\nu^*,\log G\rangle\leq \langle\pi\nu^*,
\log G^-\rangle +\delta/(2p)$. This can be seen by choosing $N$ so large 
that $\langle \nu,\log G\rangle-\langle\nu,\log G_N\rangle<\delta/(4p)$ 
and using uniform continuity of $\log G_N$ on $U^2$ to split the domain 
of integration into partition sets on which the variation of $\log G_N$ 
is less than $\delta/(4p)$.
Using this and \eqref{HIesti} in \eqref{ochapprmin}, we arrive at
$$\begin{aligned}
-\och(\phi,\lambda)&\geq \sum_{i=1}^n \lambda_i H\Bigl(\pi\mu_i^*\,\Big|\,
\int_{U_\cdot}\varphi_i\Bigr)+p\bigl(\widetilde I^2_{\pi\overline\nu^*}
(\pi\nu^*)-\langle\pi\nu^*,\log G^-\rangle\bigr)-\delta\\
&\geq \inf_{u_1,\dots,u_n\in\skrim_1(\Sigma_r)}\Bigl[\sum_{i=1}^n \lambda_i 
H\Bigl(u\,\Big|\,\int_{U_\cdot}\varphi_i\Bigr)+p\skrig^-\Bigl(\sum_{i=1}^n
\lambda_i u_i\Bigr)\Bigr]-\delta\\
&=-\widetilde\och^-(\phi,\lambda)-\delta,
\end{aligned}$$
which finishes the proof of (i).

Now we prove (ii). We choose the partition so fine that
\begin{equation}\label{partchoicea}
\bigl| \log G_M^+(\pi x,\pi y)-\log G_M(x,y)\big|<\frac\delta{2p}, \qquad 
\mbox{ for all } x,y\in U,
\end{equation}
and
\begin{equation}\label{partchoiceb}
\bigl|\log \widetilde \varphi_i(x)-\log \varphi_i(x)\bigr|<\frac\delta2,
\qquad \mbox{ for all }i\in\{1,\dots,n\}, \, x\in\supp(\varphi_i),
\end{equation}
where $\widetilde\varphi_i(x)=\sum_{l=1}^r\Leb(U_l)^{-1}
\1_{U_l}(x)\, \int_{U_l}\varphi_i$ is a discrete approximation to $\varphi_i$. 
This choice is possible since 
the functions $\log G_M$ and $\log \varphi_i$ are bounded and measurable on 
$U^2$ resp.\ on $\supp(\varphi_i)$. To every $\mu\in\skrim_1(U)$ we associate 
a $\widetilde\mu\in\skrim_1(U)$with constant density on the partition sets 
and $\pi\mu=\pi\widetilde \mu$. Note that $H(\pi\mu\mid\int_{U_\cdot}\varphi_i)
=H(\widetilde\mu\mid\widetilde \varphi_i)$ for any $i\in\{1,\dots,n\}$ and any 
$\mu\in\skrim_1(U)$. Furthermore, we easily derive from \eqref{partchoicea} 
resp.~from \eqref{partchoiceb} that 
$$
\skrig_M^+(\pi\mu)\geq \skrig(\widetilde\mu)-\frac \delta {2p}\qquad\mbox{and}
\qquad H\bigl(\pi\mu\,\big|\,\widetilde\varphi_i\bigr)\geq H(\widetilde\mu\mid\varphi_i)-\frac\delta2,
$$ 
for any $\mu\in\skrim_1(U)$ and any $i\in\{1,\dots,n\}$. Hence, we obtain, 
for any $\lambda\in{\mathbb S}_n$,
$$
\begin{aligned}
\och_M^+(\phi,\lambda)&=-\inf_{\mu_1,\dots,\mu_n\in\skrim_1(U)}
\Bigl(p\skrig_M^+\Bigl(\pi\Bigl(\sum_{i=1}^n\lambda_i\mu_i\Bigr)\Bigr)+
\sum_{i=1}^n\lambda_iH\Bigl(\pi\mu_i\Big|\int_{U_\cdot}\varphi_i\Bigr)\Bigr)\\
&\leq -\inf_{\mu_1,\dots,\mu_n\in\skrim_1(U)}\Bigl(p\skrig\Bigl(\sum_{i=1}^n
\lambda_i\widetilde\mu_i\Bigr)+\sum_{i=1}^n \lambda_i H\bigl(\widetilde\mu_i
\mid\varphi_i\bigr)\Bigr)+\delta\\
&\leq -\inf_{\mu_1,\dots,\mu_n\in\skrim_1(U)}\Bigl(p\skrig
\Bigl(\sum_{i=1}^n\lambda_i\mu_i\Bigr)+\sum_{i=1}^n \lambda_i H(\mu_i\mid\varphi_i)\Bigr)+\delta\\
&=\och(\phi,\lambda)+\delta.
\end{aligned}
$$
Finally we prove~(iii). We choose the partition so fine that~\eqref{partchoicea}
holds and such that $\mu\in H$, $\pi\mu=\pi\tmu$ imply $\tmu\in H[\delta]$.
As in the proof of~(ii) we associate to any $\mu\in H$ a measure 
$\tmu\in\skrim_1(U)$ with constant density on the partition sets and
$\pi\mu=\pi\tmu$. In particular, this implies $\tmu\in H[\delta]$.
Since $\widetilde{I}(\pi\mu)=I(\tmu)$, the statement follows as above.
\end{Proof}
\qed

\vspace{0.3cm}
We now complete the proof of Proposition~\ref{mixedmomasy}. For part~(i),
it suffices to combine \eqref{momentformula}, \eqref{limminus} and 
Lemma~\ref{fineness}(i). For the proof of part~(ii) let $\delta>0$ be small. 
By \eqref{momentformula} and Lemma~\ref{cutoff}(i) we can 
choose $M>0$ such that, for all $\lambda\in{\mathbb S}_n$,
$$
\E\Bigl[ \frac{1}{k!^{p}}\prod_{i=1}^n\langle\varphi_i,\ell\rangle^{k\lambda_i}
\Bigr] \leq \delta^k+ (1+\delta)^k \, \sum_{m=\lceil k(1-\delta)\rceil}^k 
\int_{B^{m}}dy\, \big(\Phi_{m,M}(y)\big)^p\prod_{i=1}^{n}
\prod_{j=1}^{\lambda_im} \varphi_i(y^i_j).
$$
The right hand side can further be estimated, using \eqref{limplus} and 
Lemma~\ref{fineness}(ii), for sufficiently large $k$ and all 
$\lambda\in{\mathbb S}_n$, by
$$\delta^k+ (1+\delta)^k \, \sum_{m=\lceil k(1-\delta)\rceil}^k 
e^{m(\och(\phi,\lambda)+2\delta)}.$$
Now we argue that $\och(\phi,\lambda)$ is bounded from below in
$\lambda\in{\mathbb S}_n$. Indeed, in \eqref{ochdef}, we get a
lower bound by choosing $\mu_i(dx)=c_i \phi^{2p}(x) \, dx$, and
noting that $\skrig$ is bounded from above on $\skrim_1(U)$.
{F}rom this the proof of Proposition~\ref{mixedmomasy}(2) readily 
follows.

For the proof of Proposition~\ref{mixedmomasy}(iii) let $\delta>0$ be small. 
By \eqref{genmomentformula} and Lemma~\ref{cutoff}(ii) we can 
choose $M>0$ such that, 
$$
\E\Big[ \frac{1}{k!^{p}} \ell(U)^k \1_{\{ L\in H\}} \Big] 
\leq \delta^k+ (1+\delta)^k \, \sum_{m=\lceil k(1-\delta)\rceil}^k 
\int_{U^{m}}dy\, \big(\Phi_{m,M}(y)\big)^p
\1_{\{L_{y,m}\in H[\delta]\}}
$$
Now we use \eqref{Phidisc} and note that $L_{\pi y,m}=\pi L_{y,m}$
to obtain
$$\limsup_{m\uparrow \infty}\frac 1m\log\Bigl(\int_{U^{m}}dy\,
 \big(\Phi_{m,M}(y)\big)^p \1_{\{L_{y,m}\in H[\delta]\}} \Bigr)
\leq -\inf_{\mu\in {H}[\delta]} 
\big\{ \widetilde{I}(\pi\mu)+p \skrig_M^+(\pi\mu) \big\}.$$
{F}rom here one can finish the proof of Proposition~\ref{mixedmomasy}(iii)
by an application of Lemma~\ref{fineness}(iii).

\section{Identification of the variational formula}\label{sec-anaform}

In this section we prove Proposition~\ref{propWident}. This is done in 
two steps. In Section~\ref{sec-enermeas} we identify $W(\phi)$ in terms of a 
variational problem involving energies of measures with respect 
to the Green operator on $B$. In Section~\ref{sec-enerfunc} this formula
is related to the variational formula \eqref{Theta} for $\Theta(\phi)$ in
Section~\ref{sec-enerfunc}, and this completes the proof of 
Proposition~\ref{propWident}.

\subsection{Identification of $\boldsymbol{W(\phi)}$ in terms of 
energies of measures}\label{sec-enermeas}

Recall the definition of the Green function $G$ from Section~\ref{first}
and define the associated operator $\A$ by
$$\A g(x) = \int G(x,y) g(y) \, dy, \qquad \mbox{ and } \qquad
\A \mu(x) = \int G(x,y) \, \mu(dy).$$ 
We introduce the {\em energy\/} of a measure $\mu$ on~$B$,
\begin{equation}\label{energy}
\|\mu\|_\En^2=\langle \mu,\A\mu\rangle
=\int_B\int_B\mu(\d x)G(x,y)\mu(dy),
\end{equation}
and we write $\|g\|_E=\|\mu\|_E$ if $\mu=g\, dx$.

Let $\phi=(\phi_1,\dots,\phi_n)$ be a family of nonnegative, bounded
measurable functions on $B$ having compact supports. 
The main object of this section is the variational formula
\begin{equation}\label{rho*def}
\begin{aligned}
\rho(\phi)&=\rho(\phi_1,\dots,\phi_n)\\
&=\sup\Big\{ \Big\| \sum_{i=1}^n 
\sqrt{\lambda_i} \, g_i^{2p-1} \, \phi_i \Big\|_\En^2 \, : \, 
\lambda\in\mathbb{S}_n,  g_i\in L^{2p}(B), \|g_i\|_{2p}=1\,\mbox{ for } 
i=1,\dots,n\Big\}.
\end{aligned}
\end{equation}

We first show that maximisers exist for this variational problem, and we derive the
variational equations.

\begin{lemma}[Analysis of $\rho(\phi)$]\label{solexist}
Let $\phi=(\phi_1,\dots,\phi_n)$ be a family of nonnegative, bounded
measurable functions on $B$ with compact supports.  
Then there exist $\lambda\in\mathbb{S}_n$ and $g_1,\dots,g_n\in L^{2p}(B)$ 
with $\|g_i\|_{2p}=1$ such that
\begin{equation}
\rho(\phi)=\Big\| \sum_{i=1}^n 
\sqrt{\lambda_i} \, g_i^{2p-1} \, \phi_i \Big\|_\En^2,
\end{equation}
and
\begin{equation}\label{ELenerg}
\sqrt{\lambda_i}\, \rho(\phi) \,g_i =
\phi_i \, \A\Big( \sum_{j=1}^n \sqrt{\lambda_j} 
\, g_j^{2p-1} \, \phi_j\Big),\qquad \mbox{for all }i=1,\dots,n.
\end{equation}
\end{lemma}

\begin{Proof}{Proof.}
We may assume that the supports of the $\phi_i$ are not empty.
Then it is clear that in \eqref{rho*def} we may add the conditions 
$g_i\ge 0$ and $\supp(g_i)\subset U$ for all $i=1,\dots,n$, 
where $U=\bigcup_{i=1}^n\supp(\phi_i)$ denotes the 
union of the supports of $\phi_1,\dots,\phi_n$, which is a compact subset of $B$. 
Furthermore,  we may relax the condition $\|g_i\|_{2p}=1$ to the 
condition $\|g_i\|_{2p}\leq 1$.
It is convenient to substitute $f_i=g_i^{2p-1}$ and to rewrite \eqref{rho*def} as
\begin{equation}\label{rhorewrite}
\rho(\phi)=\sup\Big\{ \Big\| \sum_{i=1}^n 
\sqrt{\lambda_i} \, f_i \, \phi_i \Big\|_\En^2 \, : \, 
\lambda\in\mathbb{S}_n, \, f_1, \ldots, f_n\in K_1 \Big\},
\end{equation}
where 
$$K_M=\{f\in L^1(U)\colon f\geq 0,\|f\|\leq M\},\qquad \mbox{ for }M>0,$$
and $\|\cdot\|=\|\cdot\|_{2p/(2p-1)}$. 
As a first step, we argue that maximisers exist for the problem 
in \eqref{rhorewrite}. In the proof of \cite[Lemma~4.3]{KM02} we showed that 
$K_1$ is weakly compact in $L^1(U)$ and that the map
$f\mapsto \|f\|_{\rm E}^2$ is upper semicontinuous on $K_1$ in the weak topology 
on $L^1(U)$. Certainly, these two statements also hold for $K_M$ for any $M>0$.
Since also the set ${\mathbb S}_n\times K_1^n$ is compact and since the map
$$
{\mathbb S}_n\times K_1^n\ni\bigl(\lambda,f_1,\dots,f_n\bigr)\mapsto \sum_{i=1}^n 
\sqrt{\lambda_i} \, f_i \, \phi_i \in K_M,
$$
(here $M>0$ is suitably chosen, only dependent on $\phi_1,\dots,\phi_n$) 
is continuous in the product topology, the existence of maximisers 
 in \eqref{rhorewrite} follows. We denote them by 
$\lambda\in{\mathbb S}_n$ and $f_1,\dots,f_n\in K_1$. It is clear that 
$\|f_i\|=1$ and $\supp(f_i)\subset\supp(\phi_i)$ 
for all $i=1,\dots,n$.

The second step is to show that $\lambda_i>0$ and that $f_i$ is bounded
away from 0 on any set where $\phi_i$ is bounded from zero, 
for any $i=1,\dots,n$. Let us first prove the first of these two statements. 
Assume the contrary, i.e., $\lambda_1=0$, say. Then we may assume that 
$f_1\phi_1$ is not almost everywhere equal to zero.
There is an $i\in\{2,\dots,n\}$ such that $\sqrt{\lambda_i}f_i\phi_i$ is not
almost everywhere equal to zero. For definiteness, we assume that $\lambda_2>0$ and that $f_2\phi_2$ is not trivial. 
With some $\delta>0$, we define
$\tilde\lambda\in\mathbb{S}_n$ by
$$
\tilde\lambda_j=\begin{cases}\delta,&\mbox{if }j=1,\\
\lambda_2-\delta,&\mbox{if }j=2,\\
\lambda_j,&\mbox{otherwise.}
\end{cases}
$$
The idea is to pick $\delta>0$ so small that 
$\|h_{\tilde\lambda}\|_{\rm E}^2>\|h_\lambda\|_{\rm E}^2$, where $h_\lambda=
\sum_{j=1}^n \sqrt{\lambda_j} \, f_j \, \phi_j$. 
This would contradict the maximality of $\lambda$ 
and therefore  prove the first assertion. We calculate
\begin{equation}\label{lambdacalc}
\begin{aligned}
\|h_{\tilde\lambda} & \|_{\rm E}^2-\|h_\lambda\|_{\rm E}^2=
\sum_{i,j=1}^n \bigl(\sqrt{\tilde\lambda_i\tilde\lambda_j}-\sqrt{\lambda_i\lambda_j}
\bigr)\bigl\langle f_i\phi_i,\A(f_j\phi_j)\bigr\rangle\\
&=\delta \bigl\langle f_1\phi_1,\A(f_1\phi_1)\bigr\rangle+
2\sqrt\delta\sqrt{\lambda_2-\delta}\bigl\langle f_1\phi_1,\A(f_2\phi_2)\bigr\rangle
+2\sqrt\delta\sum_{j=3}^n\sqrt{\lambda_j}\bigl\langle f_1\phi_1,\A(f_j\phi_j)\bigr\rangle\\
&\qquad
+2\big(\sqrt{\lambda_2-\delta}-\sqrt{\lambda_2}\big)
\sum_{j=3}^n\sqrt{\lambda_j}\bigl\langle f_2\phi_2,\A(f_j\phi_j)\bigr\rangle
 -\delta\bigl\langle f_2\phi_2,\A(f_2\phi_2)\bigr\rangle\\
&\geq 
\sqrt\delta\big(c_1\bigl\langle f_1\phi_1,\A(f_2\phi_2)\bigr\rangle-c_2\sqrt\delta),
\end{aligned}
\end{equation}
for positive constants $c_1, c_2$, not depending on $\delta$. 
Since $f_1\phi_1$ and $f_2\phi_2$ are nonnegative and not trivial, 
and since $G$ is bounded away from zero on $U^2$, it is clear that 
the right hand side of \eqref{lambdacalc} is positive for sufficiently 
small $\delta>0$. This contradicts the maximality of $\lambda$. 
Hence, $\lambda_i>0$ for all $i\in\{1,\dots,n\}$.

Now we fix a small $\delta>0$ and prove that every $f_i$ is essentially bounded
away from 0 on $\{\phi_i>\delta\}$. Abbreviate $U_1=\{\phi_1>\delta\}$ and 
assume for contradiction that $|\{f_1\leq \eps\}\cap U_1|>0$ for all
$\eps>0$. Pick some $c>0$ such that $|\{f_1>c\}\cap U_1|>0$. 
With some $a,b>0$, we define $\tilde f_1\colon U\to[0,\infty)$ by
$$\tilde f_1(x)=\begin{cases}f_1(x)+a,&\mbox{if }f_1(x)\leq\eps,\\
f_1(x)-b,&\mbox{if }f_1(x)\geq c,\\
f_1(x),&\mbox{otherwise.}
\end{cases}$$
The idea is to pick $a,b>0$ in such a way that $\|\tilde f_1\|=1$ but 
$\|\tilde h_{\lambda}\|_{\rm E}^2>\|h_\lambda\|_{\rm E}^2$, where $h_\lambda=
\sum_{j=1}^n \sqrt{\lambda_j} \, f_j \, \phi_j$, and $\tilde h_{\lambda}$ is
defined analogously with $f_1$ replaced by $\tilde f_1$. This would contradict 
the maximality of $f_1,\dots,f_n$ and therefore prove the assertion.

For notational convenience, we put $\tilde f_i=f_i$ for $i\geq 2$.
Abbreviate $\eta=1/(2p-1)$. For every sufficiently small $a$ and $\eps>0$, 
we can find $b\in(0,c/2)$ such that $\|\tilde f_1\|=1$. This implies
$$\begin{aligned}
0&=\|\tilde f_1\|^{1+\eta}-\|f_1\|^{1+\eta}\\
&=\int_{\{f\leq \eps\}\cap U_1}
\bigl[(f_1(x)+a)^{1+\eta}-f_1(x)^{1+\eta}\bigr]\,dx+
\int_{\{f\geq c\}\cap U_1}\bigl[(f_1(x)-b)^{1+\eta}-f_1(x)^{1+\eta}\bigr]\,dx.
\end{aligned}$$
Hence, for some constant $C>0$ depending neither on $a$ nor on $\eps$,
we have $b\leq Ca(a+\eps)^\eta\bigl|U_1\cap \{f_1\leq \eps\}\bigr|.$ Now 
we calculate
\begin{equation}\label{fcalc}
\begin{aligned}
\|&h_{\tilde\lambda}\|_{\rm E}^2-\|h_\lambda\|_{\rm E}^2
=\sum_{i,j=1}^n \sqrt{\lambda_i\lambda_j}\Bigl(\bigl\langle \tilde f_i\phi_i,
\A(\tilde f_j\phi_j)\bigr\rangle-\bigl\langle f_i\phi_i,\A(f_j\phi_j)
\bigr\rangle\Bigr)\\
&=\lambda_1\int_{U_1\cap\{f_1\leq \eps\}}\int_{U_1\cap\{f_1\leq \eps\}}dx\,dy\,
 G(x,y)\phi_1(x)\phi_1(y)\bigl[(f_1(x)+a)(f_1(y)+a)-f_1(x)f_1(y)\bigr]\\
&\quad+\lambda_1\int_{U_1\cap\{f_1\geq c\}}\int_{U_1\cap\{f_1\geq c\}}dx\,dy\, 
G(x,y)\phi_1(x)\phi_1(y)\bigl[(f_1(x)-b)(f_1(y)-b)-f_1(x)f_1(y)\bigr]\\
&\quad+2\lambda_1\int_{U_1\cap\{f_1\leq \eps\}}\int_{U_1\cap\{f_1\geq c\}}dx\,dy\,
 G(x,y)\phi_1(x)\phi_1(y)\bigl[(f_1(x)+a)(f_1(y)-b)-f_1(x)f_1(y)\bigr]\\
&\quad+2\sqrt{\lambda_1}\sum_{j=2}^n\sqrt{\lambda_j}\int_{U_1\cap\{f_1\leq \eps\}}
\int_{U_j}dx\,dy\, G(x,y)\phi_1(x)\phi_j(y)\bigl[(f_1(x)+a)-f_1(x)\bigr]f_j(y)\\
&\quad-2\sqrt{\lambda_1}\sum_{j=2}^n\sqrt{\lambda_j}\int_{U_1\cap\{f_1\geq c\}}
\int_{U_j}dx\,dy\, G(x,y)\phi_1(x)\phi_j(y)
\bigl[(f_1(x)-b)-f_1(x)\bigr]f_j(y)\\
&\geq aC_1|U_1\cap\{f_1\leq \eps\}|-bC_2,
\end{aligned}
\end{equation}
for some constants $C_1>0$ and $C_2>0$, neither depending 
on $a$ nor on $b$. {F}rom the bound on $b$, we see that the right hand side 
of \eqref{fcalc} is positive 
for $a>0$ and $b>0$ sufficiently small, if $\eps>0$ is sufficiently small. This 
contradicts the maximality of $f_1,\dots,f_n$. Hence, every $f_i$ is 
essentially bounded away from zero on sets of the form $\{\phi_i>\delta\}$.

The third and last step is a standard application of variational 
techniques to derive the variational equation in \eqref{ELenerg} 
for the maximisers $\lambda\in{\mathbb S}_n$ and $f_1,\dots,f_n\in K$. 
It is convenient to substitute 
$$r_i=g_i^{2p}=f_i^{2p/(2p-1)} \mbox{ for $i=1,\dots,n$, }$$
then $r_i$ is normalized in $L^1(U)$-sense. For any family of test 
functions $\varphi_i\colon\{\phi_i>\delta\}\to\R$ satisfying $\int \varphi_i=0$ 
for $i=1,\dots,n$, and for any vector $v=(v_1,\dots,v_n)$ satisfying 
$\sum_{i=1}^n v_i =0$, the objects $\lambda+\eps v$ and $r_i+\eps\varphi_i$ 
are admissible for all $\eps$ with $|\eps|$ sufficiently small, and we obtain
\begin{align}
0&=\frac{d}{d\eps}\Big|_{\eps=0} \Big\| \sum_{i=1}^n 
\sqrt{\lambda_i+\eps v_i} \, (r_i+\eps \varphi_i)^{\frac {2p-1}{2p}} 
\, \phi_i \Big\|_\En^2\nonumber\\
&=\frac {2p-1}p\sum_{i=1}^n \sqrt{\lambda_i}
\Bigl\langle\varphi_i,r_i^{-\frac 1{2p}}\phi_i,
\A\Bigl(\sum_{j=1}^n \sqrt{\lambda_j}r_j^{\frac {2p-1}{2p}}\phi_j\Bigr)\Bigr\rangle
+\sum_{i=1}^nv_i\frac 1{\sqrt{\lambda_i}}\Bigl\langle r_i^{\frac {2p-1}{2p}}\phi_i,
\A\Bigl(\sum_{j=1}^n \sqrt{\lambda_j}r_j^{\frac {2p-1}{2p}}\phi_j\Bigr)\Bigr\rangle
\nonumber\\
&=\frac {2p-1}p\sum_{i=1}^n \sqrt{\lambda_i}\bigl\langle\varphi_i, g_i^{-1}\phi_i\A(h_\lambda)\bigr\rangle
+\sum_{i=1}^nv_i\frac 1{\sqrt{\lambda_i}}\bigl\langle g_i^{2p-1}\phi_i,\A(h_\lambda)\bigr\rangle,\label{varia}
\end{align}
where we put $h_\lambda=\sum_{j=1}^n \sqrt{\lambda_j}g_j^{2p-1}\phi_j$. 
Putting $\varphi_i=0$ for all $i$, we obtain $C>0$ such that 
\begin{equation}\label{dipl}
C\sqrt{\lambda_i}=\langle g_i^{2p-1}\phi_i,\A(h_\lambda)\rangle
\mbox{ for all $i$. }\end{equation}
Multiplying this with $\sqrt{\lambda_i}$, summing over $i$ and 
using $\rho(\phi)=\langle h_\lambda,\A(h_\lambda)\rangle$,
it follows that $C=\rho(\phi)$.
Putting $v_i=0$ for all $i$ in \eqref{varia} and choosing all but one 
$\varphi_j$ equal to zero, 
we obtain the existence of $C_1,\dots,C_n>0$ such that $C_ig_i= 
\phi_i\A(h_\lambda)$ for all $i$.  Multiplying the latter
equality by $g_i^{2p-1}$, integrating over $B$ and using \eqref{dipl}, 
one easily obtains that $C_i=\rho(\phi)\sqrt{\lambda_i}$ for all $i$. 
This completes the proof of \eqref{ELenerg}.
\end{Proof}
\qed

Now we characterise $W(\phi)$ in terms of $\rho(\phi)$. Recall the 
definitions \eqref{Wdef} and \eqref{rho*def} of $W(\phi)$ 
and $\rho(\phi)$, respectively.

\begin{prop}[Relation between $W$ and $\rho$]\label{prop2}
Let $\phi=(\phi_1,\dots,\phi_n)$ be a family of nonnegative, bounded
measurable functions on $B$ having compact supports. Then
$W(\phi)=p \log \rho(\phi)$, i.e.,
\begin{equation}\label{Wrhoident}
\begin{aligned}
-\min_{\lambda\in \mathbb{S}_n} &\,\min_{\mu_1,\dots,\mu_n\in\skrim_1(B)} \, 
\bigg\{  \sum_{i=1}^n \Big\langle\lambda_i\mu_i, \log\Bigl( 
\frac{\lambda_i^p}{\phi_i^{2p}}\frac{d\mu_i}{dx}\Bigr) \Big\rangle + 
p \skrig\Big(\sum_{i=1}^n \lambda_i\mu_i \Big)\bigg\}\\
& = p \log \max\Big\{ \Big\| \sum_{i=1}^n 
\sqrt{\lambda_i} \, g_i^{2p-1} \, \phi_i \Big\|_\En^2 \, : \, 
\lambda\in\mathbb{S}_n,  g_i\in L^{2p}(B),  \|g_i\|_{2p}=1\,\mbox{ for } 
i=1,\dots,n\Big\}.
\end{aligned}
\end{equation}
An explicit one-to-one correspondence between the
maximisers on the right and the minimisers on the left hand side is given by the relation 
$g_i^{2p}=\frac{d\mu_i}{dx}$ for $i=1,\dots,n$.
\end{prop}

\begin{Proof}{Proof.}
In order to prove \lq$\leq$\rq\ in \eqref{Wrhoident}, we shall show that, for any $\lambda\in\mathbb{S}_n$
and any $\mu_1,\dots,\mu_n\in\skrim_1(B)$,
\begin{equation}\label{Wleqrho}
-\sum_{i=1}^n \Big\langle\lambda_i\mu_i, \log\Bigl( 
\frac{\lambda_i^p}{\phi_i^{2p}}\frac{d\mu_i}{dx}\Bigr) \Big\rangle -
p \skrig\Big(\sum_{i=1}^n \lambda_i\mu_i \Big)\leq p\log\|h_\lambda\|_{\rm E}^2,
\end{equation}
where we put $g_i^{2p}=\frac{d\mu_i}{dx}$ for $i=1,\dots,n$ and abbreviated $h_\lambda=\sum_{i=1}^n 
\sqrt{\lambda_i} \, g_i^{2p-1} \, \phi_i$.

Abbreviate $\mu=\sum_{i=1}^n\lambda_i\mu_i$ and $g^{2p}=d\mu/dx$. 
Using the definition of~$\skrig$, Jensen's inequality and the 
concavity of $\log$, we get the following upper bound,
\begin{align*}
-\sum_{i=1}^n & \Big\langle\lambda_i\mu_i, \log\Bigl( 
\frac{\lambda_i^p}{\phi_i^{2p}}\frac{d \mu}{dx}\Bigr) \Big\rangle -
p \skrig\Big(\sum_{i=1}^n \lambda_i\mu_i \Big) \\
& = -p \bigg[ \sum_{i=1}^n \lambda_i \Big\langle\mu_i, \log 
\frac{\lambda_i g_i^2}{\phi_i^2}\Big\rangle +
\Big\langle\mu, \log \frac{h_\lambda}{g^{2p}}\Big\rangle
- \sup_{\nu\colon \bar{\nu}=\mu}
\int \mu(dx) \int \frac{\nu(dx dy)}{\mu(dx)}
\log \frac{g^{2p}(x) G(x,y) h_\lambda(y)}{\frac{\nu(dx dy)}{dx dy}} 
\bigg]\\
& \le -p \bigg[ \sum_{i=1}^n \lambda_i \Big\langle\mu_i, \log 
\frac{\lambda_i g_i^2 h_\lambda}{\phi_i^2 g^{2p}}\Big\rangle - 
\big\langle \mu, \log \A(h_\lambda) \big\rangle \bigg]\\
& = p \bigg[ \sum_{i=1}^n \lambda_i \Big\langle g_i^{2p}, \log 
\frac{\phi_i g^{2p}}{\sqrt{\lambda_i} g_i h_\lambda}\Big\rangle +
\sum_{i=1}^n \lambda_i \Big\langle g_i^{2p}, \log \frac{\phi_i \A(h_\lambda)}
{\sqrt{\lambda_i} g_i} \Big\rangle \bigg]\\
& \le p \bigg[ \sum_{i=1}^n \lambda_i \log \Big\langle g_i^{2p-1}
\frac{\phi_i}{\sqrt{\lambda_i}},  \frac{g^{2p}}{h_\lambda}\Big\rangle +
\sum_{i=1}^n \lambda_i \log \Big\langle g_i^{2p-1}\frac{\phi_i}{\sqrt{\lambda_i}}
, {\A(h_\lambda)} \big\rangle \bigg]\\
& \le p \bigg[  \log \Big\langle \sum_{i=1}^n \sqrt{\lambda_i}\,
g_i^{2p-1}\, \phi_i, \frac{g^{2p}}{h_\lambda} \Big\rangle
+ \log\big\langle h_\lambda , \A(h_\lambda) \big\rangle
\bigg]\\
& = p \log \| h_\lambda\|^2_\En.
\end{align*}
This shows that \eqref{Wleqrho} holds and implies the upper bound in \eqref{Wrhoident}.

To prove the lower bound in \eqref{Wrhoident} we
pick, in accordance with Lemma~\ref{solexist}, maximisers 
$\lambda\in\mathbb{S}_n$ and  $g_1,\dots,g_n\in L^{2p}(B)$ for the problem of the right hand side
and show that the value of the functional of the left hand side 
for the choice  $\mu_i(dx)=g_i^{2p}(x)\,dx$ for $i=1,\dots,n$ is not
smaller than the value of the maximum on the right.
Recall from \eqref{ELenerg} that $\mu_i$-almost everywhere $\phi_i>0$.
We first find an upper bound for $\skrig(\sum_{i=1}^n\lambda_i\mu_i)$ by
picking some particular $\nu\in\skrim_1(U^2)$; recall the definition~\eqref{Gdef}
of~$\skrig$.
Indeed, define $\nu^*\in\skrim_1(U^2)$ by
$$
\nu^*(dx dy)=\frac{1}{\rho}
h_\lambda(x) \,  G(x,y) \, h_\lambda(y) \,  dx\, dy,
$$
where we abbreviated $\rho=\rho(\phi)$ and $h_\lambda=\sum_{i=1}^n 
\sqrt{\lambda_i} \, g_i^{2p-1} \, \phi_i$.
Note from \eqref{ELenerg} that, for all $i =1,\dots, n$,
\begin{equation}\label{ELen2}
\A\big(h_\lambda\big)(x)=\rho \, \sqrt{\lambda_i}\, \frac{g_i(x)}{\phi_i(x)},
\mbox{ if } \phi_i(x)>0.
\end{equation}
Note that $\nu^*$ is an admissible choice in the optimisation problem in the
definition \eqref{Gdef} of $\skrig$, by symmetry and because, using \eqref{ELen2},
$$
\overline{\nu}^*(dy)  = \frac{1}{\rho}\,
h_\lambda(y) \, \A\big(h_\lambda\big)(y)\, dy = 
\sum_{i=1}^n \lambda_i\, g_i^{2p}\, dy
= \mu(dy).
$$
Replacing the supremum over all $\nu$ by the value 
for $\nu^*$ gives the following lower bound for the left hand side 
of \eqref{Wrhoident},
\begin{align*}
-\sum_{i=1}^n & \Big\langle\lambda_i\mu_i, \log\Bigl( 
\frac{\lambda_i^p}{\phi_i^{2p}}\frac{d \mu}{dx}\Bigr)\Big\rangle - 
p \skrig\Big(\sum_{i=1}^n \lambda_i\mu_i \Big) \\
& = p \bigg[- \sum_{i=1}^n \lambda_i \Big\langle\mu_i, \log 
\frac{\lambda_i g_i^2}{\phi_i^2}\Big\rangle 
+ \sup_{\nu: \bar{\nu}=\mu} \iint \nu(dx dy) \log \frac{g^{2p}(x) G(x,y) 
g^{2p}(y)}{\frac{\nu(dx dy)}{dx dy}} \bigg]\\
& \ge p \bigg[ -\sum_{i=1}^n \lambda_i \Big\langle\mu_i, \log 
\frac{\lambda_i g_i^2}{\phi_i^2}\Big\rangle+ \iint \nu^*(dx dy)
\log \frac{g^{2p}(x) G(x,y) g^{2p}(y)}{\frac{\nu^*(dx dy)}{dx dy}} 
\bigg]\\
& = p \bigg[- 2\,\sum_{i=1}^n \lambda_i \Big\langle\mu_i, \log 
\frac{\sqrt{\lambda_i} g_i}{\phi_i}\Big\rangle +\log \rho\,  
+ 2 \int \mu(dx)\log \frac{g^{2p}(x)}{h_\lambda(x)} \bigg]\\
& = p \log \|h_\lambda\|_{\rm E}^2= p \log \rho,
\end{align*}
because $\sqrt{\lambda_i} g_i h_\lambda= g^{2p} \phi_i$ by \eqref{ELen2}
and the definition of $h_\lambda$. This completes the proof.
\end{Proof}
\qed

\subsection{Identification of $\boldsymbol{W(\phi)}$ in 
terms of energies of functions}\label{sec-enerfunc}

In this section, we identify the variational formula \eqref{rho*def} for 
$\rho(\phi)$ in terms of the formula \eqref{Theta} for 
$\Theta(\phi)$ and prove Proposition~\ref{analtheta}. 
As a first step, we prove that
minimisers exist in \eqref{Theta}, and we derive their
variational equation.

\begin{lemma}[Analysis of $\Theta(\phi)$]\label{solexauch}
Let $\phi=(\phi_1,\dots,\phi_n)$ be a family of nonnegative, bounded
measurable functions on $B$ having compact supports. Then there
exists a $\psi\in\skrid(B)$, which satisfies
$$
\Theta(\phi)=\frac p2 \|\nabla \psi\|^2_2,
\mbox{ and }  \, \sum_{i=1}^n \|\phi_i \psi\|_{2p}^2 =1,
$$
and with $h=\sum_{i=1}^n \|\phi_i \psi\|_{2p}^{2-2p}\,\phi_i^{2p}$
we have the variational equations
\begin{equation}\label{EL2}
\frac{p}{\Theta(\phi)}\, \psi = \A( \psi^{2p-1} \,h )\quad\mbox{and}\quad
-\frac p2\Delta\psi={\Theta(\phi)} \psi^{2p-1}\,h.
\end{equation}
\end{lemma}

\begin{Proof}{Proof.}
As a first step, we derive the existence of a minimiser in \eqref{Theta}.
Let $(\psi_k\, : \, k\in\N)$ be a minimising sequence, that is, 
the functions $\psi_k\in \skrid(B)$ are nonnegative and satisfy
$\sum_{i=1}^n \|\phi_i \psi_k\|_{2p}^2 =1$ for any $k\in\N$, and 
$\lim_{k\uparrow\infty}\frac 12\|\nabla\psi_k\|_2^2 =\Theta(\phi)$.

Let $\psi_*\in \skrid(B)$ 
denote the weak limit of a subsequence in accordance with Lemma~\ref{subsequences}. 
By local strong convergence in $L^{2p}(B)$, we also have 
$\sum_{i=1}^n \|\phi_i \psi_*\|_{2p}^2 =1$.
By weak lower semicontinuity of the map $\psi\mapsto \|\nabla\psi\|_2^2$ 
(apply \cite[Theorem 2.11]{LL97}), we have that $\frac 12\|\nabla\psi_*\|_2^2
\leq \liminf_{k\uparrow\infty}\frac 12\|\nabla\psi_k\|_2^2=\Theta(\phi)$. 
Since  $\psi_*$ is certainly nonnegative, it is a minimiser 
in \eqref{Theta}.

The second step is the derivation of the variational equation in 
\eqref{EL2} for any minimiser $\psi_*$ in \eqref{Theta}. 
Since $\|\nabla|\psi|\|_2^2 = \|\nabla\psi\|_2^2$ for any 
$\psi\in \skrid(B)$ (see \cite[Theorem~6.17]{LL97}), and since 
$\sum_{i=1}^n\|\phi_i \psi\|_{2p}^2 $ is positive homogeneous of order two in $\psi$,
$\psi_*$ is also a minimiser in the variational problem 
\begin{equation}\label{Thetarepr}
\Theta(\phi)=\inf_{\psi\in\skrid(B)}
\frac{\frac p2 \|\nabla \psi\|^2_2}{\sum_{i=1}^n \| \phi_i \psi \|_{2p}^2}.
\end{equation}
Denote the quotient on the right hand side of \eqref{Thetarepr} by $F(\psi)$. Let 
$\varphi\in C_{\rm c}^\infty(B)$ be a smooth test function, then
the map $\eps\mapsto F(\psi_*+\eps\varphi)$ can easily be differentiated at $\eps=0$. By minimality
of $\psi_*$ for $F$, this derivative is equal to zero. Recalling that 
$\sum_{i=1}^n \|\phi_i \psi_*\|_{2p}^2=1$, this implies that
\begin{equation}\label{ELeq}
\begin{aligned}
0&=\frac{d}{d\eps}\Big|_{\eps=0}\bigl\|\nabla(\psi_*+\eps\varphi)\bigr\|_2^2
-\sum_{i=1}^n \|\nabla\psi_*\|_2^2\frac{d}{d\eps}\Big|_{\eps=0}\|\phi_i (\psi_*+\eps\varphi)\|_{2p}^2\\
&=2\int_B\nabla\psi_*\cdot\nabla\varphi-4\frac{\Theta(\phi)}p\sum_{i=1}^n\|\phi_i\psi_*\|_{2p}^{2-2p}
\bigl\langle\varphi, \phi_i^{2p}\psi_*^{2p-1}\bigr\rangle\\
&=-2\Bigl\langle \varphi,\Delta \psi_*+2\frac{\Theta(\phi)}p\psi_*^{2p-1}h\Bigr\rangle,
\end{aligned}
\end{equation}
where we used the definition of the distributional Laplacian in the last step, 
and  $h=\sum_{i=1}^n \|\phi_i\psi_*\|_{2p}^{2-2p}
\phi_i^{2p}$ as in \eqref{EL2}.
As \eqref{ELeq} holds for any smooth test function $\varphi$, we infer that 
the function in the right argument of the brackets on the right hand side is
equal to zero, i.e., $-\frac 12\Delta\psi_*=\frac{\Theta(\phi)}p \psi_*^{2p-1}h$, which is the second identity in 
\eqref{EL2}. By \cite[Th.~6.21]{LL97}, the function 
${\psi}=\frac{\Theta(\phi)}p\A(\psi_*^{2p-1} h)$ 
satifies $-\frac 12 \Delta \psi= \frac{\Theta(\phi)}p \psi_*^{2p-1}h$.
Hence, by \cite[Th.~9.3]{LL97}, ${\psi}$ differs from $\psi_*$ by a 
harmonic function in $\skrid(B)$, which therefore vanishes. 
This ends the proof of \eqref{EL2}.
\end{Proof}\qed

Now we identify $\rho(\phi)$ in terms of $\Theta(\phi)$. The following proposition
completes the proof of Theorem~\ref{momasy}, with the help of
Proposition~\ref{asym} and Proposition~\ref{prop2}.

\begin{prop}[Relation between $\rho$ and $\Theta$]\label{prop3}
Let $\phi=(\phi_1,\dots,\phi_n)$ be a family of nonnegative, bounded
measurable functions on $B$ with compact supports. Then 
$\rho(\phi)=p/\Theta(\phi)$, i.e.,
\begin{equation}\label{rhoThetaident}
\begin{aligned}
\max\Big\{ &\Big\|  \sum_{i=1}^n 
\sqrt{\lambda_i} \,  g_i^{2p-1} \, \phi_i \Big\|_\En^2 \, : \,
\lambda\in\mathbb{S}_n,  g_i\in L^{2p}(B),\,  \|g_i\|_{2p}=1\,
\mbox{ for } i=1,\dots,n\Big\}^{-1} \\
& = \min\Big\{ \sfrac 12 \|\nabla \psi\|^2_2 \, : \, \psi\in\skrid(B), \,
\sum_{i=1}^n \|\phi_i \psi\|_{2p}^2 =1 \Big\}.
\end{aligned}
\end{equation}
\end{prop}

\begin{remark}
The proof gives an explicit one-to-one correspondence between the
maximisers on the left and the minimisers on the right hand side, 
see \eqref{geq} and \eqref{leq}, respectively.\hfill{$\Diamond$}
\end{remark}

\begin{Proof}{Proof.}
For the proofs of both \lq$\geq$\rq\ and \lq$\leq$\rq\ in \eqref{rhoThetaident}, we
pick the maximiser resp.\ the minimiser in one variational formula,
construct admissible objects for the other one, and show that the other
functional attains the inverse of the value of the maximum resp.\ minimum. 

Let us begin
with the proof of \lq$\geq$\rq.
Pick maximisers $\lambda\in\mathbb{S}_n$ and $g_1,\dots,g_n\in L^{2p}(B)$ of the 
formula on the right hand side of \eqref{rhoThetaident} in accordance with 
Lemma~\ref{solexist}. Define
\begin{equation}\label{geq}
\psi= \frac 1{\rho(\phi)} \,\, \A\Big( \sum_{j=1}^n
\sqrt{\lambda_j} g_j^{2p-1} \phi_j\Big).
\end{equation}
Then, by \eqref{ELenerg}, for all $i=1,\dots,n$,
\begin{equation}\label{ELen3}
\psi(x)= \sqrt{\lambda_i} \frac{g_i}{\phi_i}(x) \mbox{ for all } 
\phi_i(x)>0.
\end{equation}
Hence,
$$
\sum_{i=1}^n \|\phi_i \psi\|_{2p}^2
= \sum_{i=1}^n \lambda_i \int g_i^{2p}(x) \, dx =1.
$$
Then the energy of the measure $\frac 1{\rho(\phi)}\sum_{j=1}^n
\sqrt{\lambda_j} g_j^{2p-1} \phi_j\, dx$ can be calculated as follows.
$$
\begin{aligned}
\Bigl\|\frac 1{\rho(\phi)} \sum_{j=1}^n
\sqrt{\lambda_j} g_j^{2p-1} \phi_j\Bigr\|_{\rm E}^2 & =\frac 1{\rho(\phi)}
\int_B dx\, \psi(x) \, \sum_{j=1}^n \sqrt{\lambda_j} g_j^{2p-1}(x)\, \phi_j(x)\\
& =  \frac 1{\rho(\phi)}\, \int_B dx\,
\, \sum_{j=1}^n \lambda_j g_j^{2p}(x) 
=  \frac 1{\rho(\phi)}.
\end{aligned}
$$
By Lemma~\ref{energies} we have that $\psi\in\skrid(B)$ and that the energy 
of $\psi$ equals the energy of the measure $\frac 1{\rho(\phi)}\sum_{j=1}^n
\sqrt{\lambda_j} g_j^{2p-1} \phi_j\, dx$, i.e.,\ $\frac 12  \|\nabla \psi\|^2_2=  
\frac 1{\rho(\phi)}$. This implies \lq$\geq$\rq. 

To prove  \lq$\leq$\rq, we choose $\psi$ as the minimiser of the problem on
the right hand side of in \eqref{rhoThetaident}, by Lemma~\ref{solexauch}. 
We define $g_1,\dots,g_n\in L^{2p}(B)$ and $\lambda\in \mathbb{S}_n$ by 
\begin{equation}\label{leq}
g_i=\frac{\psi \phi_i}{\|\psi \phi_i\|_{2p}}\qquad\mbox{and}\qquad
\lambda_i=\|\psi \phi_i\|_{2p}^2\quad\mbox{for }i=1,\dots,n.
\end{equation}
Note that $\|g_i\|_{2p}=1$ for all $i=1,\dots,n$ and that $\lambda_1,\dots,\lambda_n$ 
are nonnegative numbers summing to one. Hence, $g_1,\dots,g_n$ and $\lambda=(\lambda_1,\dots,\lambda_n)$
are admissible for the formula on the left hand side of  \eqref{rhoThetaident}.
We find, using the first identity in \eqref{EL2},
\begin{align*}
\Big\| \sum_{i=1}^n  & \sqrt{\lambda_i} g_i^{2p-1} \phi_i \Big\|^2_\En =
\Big\|\sum_{i=1}^n\|\psi\phi_i\|_{2p}^{2-2p}
\, \psi^{2p-1} \phi_i^{2p}\Big\|^2_\En \\
& = \sum_{j=1}^n \Big\langle\psi^{2p-1}\, \phi_j^{2p} \|\psi\phi_j\|_{2p}^{2-2p}\, ,
\A\Big( \psi^{2p-1} \,  \sum_{i=1}^n \phi_i^{2p}\, \|\psi\phi_i\|_{2p}^{2-2p}\Big) 
\Big\rangle \\
&= \frac p{ \Theta(\phi)}\, \sum_{j=1}^n
\Big\langle\psi^{2p-1}\phi_j^{2p}\,\|\psi\phi_j\|_{2p}^{2-2p} , \psi 
\Big\rangle=\frac{p}{ \Theta(\phi)}.
\end{align*}
This completes the proof of the proposition.\end{Proof}\qed

\section{Appendix: The space ${\skrid(B)}$}\label{space}

We now recall the definition of the function space $\skrid(B)$ 
and state some properties of this space. All of this material 
is known to the experts, but we find it convenient to collect some 
technical facts which are used at some places.

In the case of $B$ bounded, $\skrid(B)$ is the classical Sobolev space 
$H_0^1(B)$ which is defined as the closure of $\skric_{\rm c}^\infty(B)$ 
in the sense of the Sobolev norm 
$\psi\mapsto (\|\nabla\psi\|_2^2+\|\psi\|_2^2)^{1/2}$ in the Sobolev
space $H^1(B)$. We first give a relation between $H_0^1(B)$ and $H^1(\R^d)$ 
in the case of a $C^1$-boundary.

\begin{lemma}\label{HaNull} Let $B\subset \R^d$ be an open bounded set with $C^1$-boundary.
Let $\psi\in H^1(\R^d)$ such that $\psi=0$ a.e.\ on $B^{\rm c}$. Then the 
restriction of $\psi$ to $B$ lies in $H_0^1(B)$.
\end{lemma}

\begin{Proof}{Proof.}
Our proof is an adaptation of the proof of Theorem~3 in \cite[Section~5.3.3]{Ev98}.
First we pick, for any $\eps>0$, a function $\varphi_\eps\in\skric^\infty(\R^d)$ such
that $\varphi_\eps\to\psi$ as $\eps\downarrow 0$ in the Sobolev norm, and
such that $\supp(\varphi_\eps)$ is contained in the open $\eps$-neighbourhood of 
$B$, which we denote by $B[\eps]$.

Fix $x^0\in\partial B$. Since $\partial B$ is $C^1$, there are $r>0$ and a 
$C^1$-function $\gamma\colon \R^{d-1}\to\R$ such that 
$$B\cap B(x^0,r)=\{x\in B(x^0,r)\colon x_1<\gamma(x_2,\dots,x_d)\}.$$
Let $V=B\cap B(x^0,r/2)$. For $\eps>0$ define $\psi_\eps\colon V\to\R$ by 
$\psi_\eps(x)=\varphi_\eps(x+\eps\lambda{\rm e}_1)$ for $x\in V$, 
where ${\rm e}_1$ denotes the first unit vector and $\lambda>1$ is chosen
such that $\supp(\psi_\eps)\subset U$. Then $\psi_\eps\in\skric^\infty(V)$.
Now the continuity of the $L^2$-norm under translations shows that 
$\lim_{\eps\downarrow 0}\|\psi_\eps-\widetilde\psi_\eps\|_{H^1(V)}=0$.
In particular, we have that $\psi_\eps\to\psi$ in $H^1(V)$.

Now the rest of the proof is as in the proof of Theorem~3 in 
\cite[Section~5.3.3]{Ev98}.
Indeed, using the compactness of $\partial B$, we find a finite covering of 
$\partial B$ with balls $B_1,\dots,B_N$ in which $\partial B$ can be mapped 
differentiably onto a hyperplane.  Within the ball $B_i$,
we can approximate $\psi$ for any $i$ in $H^1$-norm by a $C^\infty$-function
$\psi^{\ssup{i}}$ whose support lies 
within $B$. Extend the covering of $\partial B$ to a covering of $B$ by adding
a suitable open set $B_0$ whose closure is contained in $B$. On $B_0$, we can
approximate $\psi$ in $H^1$-norm by a $C^\infty$-function 
$\psi^{\ssup{0}}$ with support within $B$ (use Theorem~1 in 
\cite[Section~5.3.3]{Ev98}). Now pick a smooth partition of the unity,
$(\zeta_i \colon  i=0,\dots,N)$, subordinated to the covering 
$B_0,\dots,B_N$ of $B$, and put $\varphi=\sum_{i=0}^N\psi^{\ssup{i}}\zeta_i$. 
It is then easily seen that $\varphi$ lies in $\skric^\infty_{\rm c}(B)$ and 
approximates $\psi$ in $H^1(B)$-norm. This completes the proof.
\end{Proof}
\qed

In the case that $B=\R^d$, the space $\skrid(\R^d)=D^1(\R^d)$ is the space 
of functions $f\in L^1_{\rm loc}(\R^d)$, which 
vanish at infinity, i.e., $\{x\in\R^d\colon |f(x)|>a\}$ has finite Lebesgue 
measure for any $a>0$, and whose distributional gradient is in $L^2(\R^d)$.
Now we collect some sequential compactness properties of the space 
$\skrid(B)$.

\begin{lemma}\label{subsequences}
Suppose $(\psi_k)_{k\in\N}$ is a sequence in $\skrid(B)$ such that 
$(\|\nabla\psi_k\|_2)_{k\in\N}$
is bounded.  Fix any $q\in(1,2d/(d-2))$ for $d\geq 3$ and any $q>1$ for $d\le 2$.
Then there exists $\psi\in\skrid(B)$ and a subsequence $(\psi_{k_j})_{j\in\N}$
such that $\nabla\psi_{k_j}\to\nabla\psi$ weakly in $L^2(B)$ and $\psi_{k_j}\to\psi$
locally strongly in $L^q(B)$.
\end{lemma}

\begin{Proof}{Proof.}
Let us recall standard Sobolev inequalities, see \cite[Theorems~8.3, 8.5]{LL97}. 
There are positive constants $S_d$ for $d\geq 3$ and $S_{2,r}$ for $r>2$ 
such that 
\begin{equation}\label{Sobolev}
\begin{array}{rcll}
S_d\|\psi \|_{2d/(d-2)}^2&\leq &\|\nabla \psi \|_2^2,\qquad&\mbox{for }d\geq 3, \psi \in D^1(\R^d),\\
S_{2,r}\|\psi \|_{r}^2&\leq& \|\nabla \psi \|_2^2+\|\psi \|_2^2,
\qquad&\mbox{for }d=2, \psi \in H^1(\R^d), r> 2.
\end{array}
\end{equation}

We first consider the case $B=\R^d$. In particular, $d\geq 3$.
Fix $1<q<2d/(d-2)$ and apply H\"older's inequality and Sobolev's inequality to get,
for any bounded measurable set $A\subset\R^d$ and any $\psi\in D^1(\R^d)$,
\begin{equation}\label{HoelderRd}
\|\psi \1_A\|_q\leq \|\psi\|_{2d/(d-2)}\Leb(A)^{\frac{2d-dq+2q}{2dq}}
\leq \|\nabla\psi\|_2 S_d^{-1/2}\Leb(A)^{\frac{2d-dq+2q}{2dq}},
\end{equation}
where $\Leb$ denotes the Lebesgue measure.

Now suppose that $(\psi_k)_{k\in\N}$ is a sequence in $D^1(\R^d)$ such that $(\|\nabla\psi_k\|_2)_{k\in\N}$
is bounded. The estimate in \eqref{HoelderRd} shows that $(\psi_k)_{k\in\N}$
is locally bounded in $L^q(\R^d)$. By the Banach-Alaoglu Theorem, there is a subsequence
$(\psi_{k_j})_{j\in\N}$ which converges to some $\psi\in L^q(\R^d)$ locally weakly in $L^q(\R^d)$
and $\nabla\psi_{k_j}\to u$ weakly in $L^2(\R^d)$ for some $u\in L^2(\R^d)$. By \cite[Theorem~8.6]{LL97}, 
the subsequence converges even locally strongly in $L^{q}(\R^d)$, and 
$u=\nabla \psi$. This completes the proof for $B=\R^d$.

Now we turn to the case of bounded $B$, for general $d\geq 1$. For
any $\psi\in L^2(B)$ we define the extension $\psi_*\in L^2(\R^d)$
by $\psi_*(x)=0$ for $x\not\in B$ and $\psi_*(x)=\psi(x)$ for $x\in B$.
Then, for $\psi\in H^1_0(B)$ we have $\nabla\psi_*=(\nabla\psi)_*$. 
Indeed, let $\varphi\in C_c^\infty(B)$ be any test function and 
$\psi_n\in C_c^\infty(B)$ be a sequence of functions approximating 
$\psi$ in the norm of $H^1_0(B)$. Then, applying partial integration 
to test functions, for all $1\le i\le d$,
\begin{align*}
\int \varphi \, \big(\frac{\partial \psi}{\partial x_j}\big)_*
 = \int_B \varphi \, \frac{\partial \psi}{\partial x_j}
 =  \lim_{n\uparrow\infty} \int_B \varphi \, \frac{\partial \psi_n}{\partial x_j}
 = \lim_{n\uparrow\infty} \int_B \, \frac{\partial\varphi}{\partial x_j} \psi_n
 = \int_B \frac{\partial\varphi}{\partial x_j} \, \psi.
 = \int \frac{\partial\varphi}{\partial x_j} \, \psi_*.
\end{align*}
This fact will be used in the sequel mostly without further 
notice. It implies, for instance, that Sobolev's inequality \eqref{Sobolev} 
is applicable to functions in $\skrid(B)=H_0^1(B)$. 

Suppose that $(\psi_k)_{k\in\N}$ is a sequence in $H_0^1(B)$ such that $(\|\nabla\psi_k\|_2)_{k\in\N}$
is bounded. In the case $d\geq 3$ and $1<q<2d/(d-2)$, similarly to \eqref{HoelderRd}, one derives that 
$(\psi_k)_{k\in\N}$ is bounded in $L^q(B)$. By the Banach-Alaoglu theorem 
in the space $H_0^1(B)$, a subsequence converges weakly to some 
$\psi\in H_0^1(B)$, and the rest of the proof is as above in the case $B=\R^d$.

In the case $d\le 2$, fixing any $q>1$, we first argue that there is a constant 
$C>0$ (depending only on $B$ and $q$) such that
\begin{equation}\label{Hoelderd=2}
\|\psi\|_q\leq C \|\nabla\psi\|_2,\qquad\mbox{for any }\psi\in H_0^1(B).
\end{equation}
In order to prove \eqref{Hoelderd=2} in $d=2$, use H\"older's inequality and the 
Sobolev inequality in~\eqref{Sobolev} to obtain, for any $\psi\in H_0^1(B)$ and 
any $r>2$,
\begin{equation}\label{d=2esti}
\|\psi\|_2^2\leq \|\psi\|^2_r\Leb(B)^{1-2/r}
\leq\frac{\Leb(B)^{1-2/r}}{S_{2,r}}\bigl( 
\|\nabla \psi \|_2^2+\|\psi \|_2^2\bigr).
\end{equation}
It is known \cite[Theorem~8.5]{LL97} that $1/S_{2,r}<(r^2 (r-2)/
[(r-1)8\pi])^{1-2/r}/(r-1)$. A Taylor approximation for $r\downarrow 2$ 
shows that the quotient on the right side of \eqref{d=2esti}
is smaller than one for $r>2$ sufficiently close to 2. For this $r$, \eqref{d=2esti} 
can be solved for $\|\psi\|_{2}^2$, and we obtain the existence of a constant $c>0$ such that 
$\|\psi\|_{2}^2\leq c\|\nabla\psi\|_2^2$. Use this estimate on the right hand side of \eqref{Sobolev} 
for $r=q$ to obtain that \eqref{Hoelderd=2} holds for some $C>0$, only depending on
$B$ and $q$. 

In order to prove \eqref{Hoelderd=2} in $d=1$, we use the simple inequality
\begin{equation}\label{primiSob}
|f(x)|^2 \le \|f'\|_2\,\|f\|_2, \qquad\mbox{ for } f\in H^1(\R), x\in \R,
\end{equation}
see e.g.~\cite[Theorem 8.5(6)]{LL97}. Now assume $q>2$. Raising \eqref{primiSob}
to the power $q/2$ and integrating over $B$, we get, for $f\in H^1_0(B)$,
$$\|f\|_q^2 \le \|f'\|_2\,\|f\|_2\,\Leb(B)^{2/q}\le
\|f'\|_2\,\|f\|_q\,\Leb(B)^{1/q+1/2},$$
where we used H\"older's inequality in the second step. This shows
\eqref{Hoelderd=2} in the case $q>2$. The general case follows by a
further application of H\"older's inequality.

The remainder of the proof in the case $d\le 2$ is as above.
\end{Proof}\qed

Recall the definition of the energy of a measure from \eqref{energy}. The following 
connection between the energy of functions in $\skrid(B)$ and the energy of 
measures will be important.

\begin{lemma}\label{energies} 
For any (positive) absolutely continuous measure $\mu\in\skrim(B)$ 
whose support is a 
compact subset of $B$ and whose energy $\|\mu\|_\En^2$ is finite, the 
function $\psi=\A (\mu)$ lies in $\skrid(B)$ and satisfies 
$\frac 12 \|\nabla\psi\|_2^2=\|\mu\|_\En^2$.
\end{lemma}

\begin{Proof}{Proof.}
Let us look at bounded $B$ first. {F}rom \eqref{HoelderRd} in the case $d\ge 3$
and from \eqref{Hoelderd=2} in the case $d=2$ we get that for some $c>0$ we have
$\|\nabla f\|_2^2 \ge c \|f\|_2^2$, for all $f\in H_{0}^{1}(B)$. By \cite[Proposition~2.5.1]{AS98}, this 
\emph{coercivity condition\/} implies that $\psi\in H_{0}^{1}(B)$ and
$\frac 12 \|\nabla\psi\|_2^2=\|\mu\|_\En^2$, as claimed.

Suppose now that $B=\R^d$, $d\ge 3$. Choose $n\in\N$ so large that 
the open centred ball $B(0,n)$ contains
the support of $\mu$. The first part shows that the function $\psi_n=\A_n (\mu)$,
defined with the operator $\A_n$ associated with the Green function $G_n$
on $B(0,n)$, lies in $H_0^1(B(0,n))$ and satisfies 
$\frac 12 \|\nabla\psi_n\|_2^2=\|\mu\|_{\En,n}^2$, 
where the energy $\|\cdot\|_{\En,n}^2$ is taken with 
respect to the domain $B(0,n)$. We can extend each $\psi_n$ by zero to the whole
of $\R^n$ and call the extension $\psi_n$ again. As $n\uparrow \infty$ we have 
$\psi_n=\A_n \mu\uparrow \A\mu=\psi$ and $\|\nabla\psi_n\|_2^2=\langle \mu,\A_n\mu\rangle
\uparrow \langle \mu,\A\mu\rangle= \|\mu\|_\En^2$.
{F}rom this, in combination with Lemma~\ref{subsequences}, we see that $\psi\in D^1(\R^d)$.

Finally, we have to show that $\lim_{n\uparrow\infty}\|\nabla\psi_n\|_2=\|\nabla\psi\|_2$. 
For this, it sufficient to show that 
$\lim_{n,m\uparrow\infty}\|\nabla(\psi_n-\psi_m)\|_2^2=0$.
Use partial integration (see \cite[Theorem~6.21]{LL97}) and the facts that 
$-\frac 12 \Delta \psi_n=\mu$ on $B(0,n)$ (as in \cite[Theorem~6.21]{LL97})
and $\psi_m=\A_m\mu$ to see that, for any $n>m$,
\begin{align*}
\sfrac 12 \|\nabla(\psi_n-\psi_m)\|
&=\sfrac 12 \|\nabla \psi_n\|_2^2+\sfrac 12 \|\nabla \psi_m\|_2^2-
\int_{B(0,n)}\nabla\psi_n\cdot\nabla\psi_m
=\|\mu\|_{\En,n}^2+\|\mu\|_{\En,m}^2+\langle \Delta\psi_n,\psi_m\rangle\\
&=\langle\mu,\A_n \mu\rangle+\langle\mu,\A_m \mu\rangle-2\langle\mu,\A_m\mu\rangle
=\langle\mu,(\A_n-\A_m)\mu\rangle\\
&=\iint\mu(\d x)\bigl(G_n(x,y)-G_m(x,y)\bigr)\mu(\d y),
\end{align*}
where $G_n$ denotes Green function of $B(0,n)$.
By Lebesgue's theorem, the right hand side vanishes as $n,m\uparrow\infty$.
\end{Proof}
\qed

{\bf Acknowledgement:} This work was supported in part by DFG grant contract 
number 234298, and by grant NAL/00631/G from the Nuffield foundation.


\bigskip

\textsc{%
\begin{tabular}
[c]{lp{0.3cm}l}%
Technische Universit\"at Berlin& & University of Bath\\
Institut f\"ur Mathematik & & Department of Mathematical Sciences\\
Strasse des 17.~Juni 136& & Claverton Down\\
10623 Berlin & & Bath BA2 7AY\\
Germany. & & United Kingdom.\\
E-Mail: {\rm koenig@math.tu-berlin.de}
& & E-Mail: {\rm maspm@bath.ac.uk } \\
\end{tabular}}

\end{document}